\documentclass[letterpaper,11pt,reqno,onefignum,onetabnum]{amsart}
\usepackage[english]{babel}
\usepackage{amsmath}
\usepackage{amsthm}
\usepackage{verbatim}
\usepackage{mathrsfs}
\usepackage{bm}
\usepackage[foot]{amsaddr}
\usepackage[dvipsnames]{xcolor}
\usepackage{hyperref}
\hypersetup{
	colorlinks = true,
	linkcolor = OliveGreen,
	anchorcolor = OliveGreen,
	citecolor = OliveGreen,
	filecolor = OliveGreen,
	urlcolor = OliveGreen
}
\usepackage{algorithm}
\usepackage{algorithmic}
\usepackage{float}
\usepackage{lipsum}
\usepackage{amsfonts}
\usepackage{amssymb}
\usepackage{graphicx}
\usepackage{epstopdf}
\usepackage{cases}
\usepackage{multirow}
\ifpdf
\DeclareGraphicsExtensions{.eps,.pdf,.png,.jpg}
\else
\DeclareGraphicsExtensions{.eps}
\fi
\usepackage{verbatim}
\usepackage{mathrsfs}
\usepackage{bm}

\newtheorem{teo}{Theorem}[section]
\newtheorem{prop}[teo]{Proposition}
\newtheorem{lem}[teo]{Lemma}

\newtheorem{pro}[teo]{Problem}
\newtheorem{algo}[teo]{Algorithm}

\newtheorem{rem}[teo]{Remark}
\newtheorem{exam}[teo]{Example}
\newtheorem{line search}{Line Search}

\usepackage{lipsum}
\usepackage{amsfonts}
\usepackage{amssymb}
\usepackage{graphicx}
\usepackage{epstopdf}

\usepackage{cases}
\usepackage{multirow}
\usepackage{tikz}
\usepackage{booktabs}%
\usepackage{hhline}
\usetikzlibrary{matrix}
\usetikzlibrary{arrows}
\ifpdf
\DeclareGraphicsExtensions{.eps,.pdf,.png,.jpg}
\else
\DeclareGraphicsExtensions{.eps}
\fi
\usepackage{verbatim}
\usepackage{mathrsfs}
\usepackage{bm}
\usepackage{color}
\usepackage{xcolor}
\usepackage{caption}
\usepackage{subcaption}
\usepackage{enumitem}

\newcommand{\N}{\mathbb N}

\newcommand{\R}{\mathbb R}

\renewcommand{\H}{\mathcal{H}}
\newcommand{\G}{\mathcal G}
\newcommand{\K}{\mathcal K}

\newcommand{\HH}{{\bm{\mathcal{H}}}}

\newcommand{\id}{\textnormal{Id}}
\newcommand{\x}{\bm x}

\newcommand{\fu}{{\textnormal{\textbf{u}}}}
\newcommand{\fz}{{\textnormal{\textbf{z}}}}
\newcommand{\weak}{\rightharpoonup}
\newcommand{\ran}{\textnormal{ran}\,}
\newcommand{\dom}{\textnormal{dom}\,}

\newcommand{\zer}{\textnormal{zer}}

\newcommand{\gra}{\textnormal{gra}\,}

\newcommand{\lb}[1]{%
  \ifmmode
    \textcolor{red}{#1}%
  \else
    {\color{red}#1}%
  \fi
}

\newcommand{\scal}[2]{{\left\langle{{#1}\mid{#2}}\right\rangle}}

\newcommand{\menge}[2]{\big\{{#1}~\big |~{#2}\big\}}
\newcommand{\pinf}{\ensuremath{{+\infty}}}
\newcommand{\minf}{\ensuremath{{-\infty}}}

\newcommand{\RP}{\ensuremath{\left[0,+\infty\right[}}

\newcommand{\RPP}{\ensuremath{\left]0,+\infty\right[}}

\newcommand{\RX}{\ensuremath{\left]-\infty,+\infty\right]}}

\newcommand{\sri}{\ensuremath{\text{\rm sri}\,}}
\newcommand{\inte}{\ensuremath{\operatorname{int}}}

\newcommand{\prox}{\ensuremath{\text{\rm prox}\,}}
\newcommand{\conv}{\ensuremath{\text{\rm conv}\,}}

\newcommand{\minimize}[2]{\ensuremath{\underset{\substack{{#1}}}%
		{\mathrm{minimize}}\;\;#2 }}

\usepackage{geometry}
\numberwithin{equation}{section}

\numberwithin{equation}{section}

\floatname{algorithm}{Line Search}

\DeclareFontEncoding{FMS}{}{}
\DeclareFontSubstitution{FMS}{futm}{m}{n}
\DeclareFontEncoding{FMX}{}{}
\DeclareFontSubstitution{FMX}{futm}{m}{n}
\DeclareSymbolFont{fouriersymbols}{FMS}{futm}{m}{n}
\DeclareSymbolFont{fourierlargesymbols}{FMX}{futm}{m}{n}
\DeclareMathDelimiter{\nr}{\mathord}{fouriersymbols}{152}{fourierlargesymbols}{147}

\DeclareMathOperator*{\argmin}{arg\,min}
\DeclareMathDelimiter{\nr}{\mathord}{fouriersymbols}{152}{fourierlargesymbols}{147}
\DeclareMathAlphabet{\mathpzc}{OT1}{pzc}{m}{it}

\title[Primal-dual Line Search for Splitting Algorithms]{\vspace*{-2.5cm}Primal-dual Line Search for Splitting Algorithms beyond Lipschitz Continuity}

\author{Felipe Atenas$^1$}
\author{Luis Briceño-Arias$^2$}
\author{Patricio Guerrero$^3$}
\author{Fernando Rold\'an$^{4}$}
\address{$^1$ Centro de Modelamiento Matem\'atico (CNRS IRL2807), Universidad de Chile, Santiago, Chile. {\it 
	E-mail address: } 
	{\sf{fatenas@cmm.uchile.cl}}.}
\address{$^2$ Departamento de Matemática, Universidad Técnica Federico Santa María,
Santiago, Chile. {\it 
	E-mail address: } 
	{\sf{luis.briceno@usm.cl}}. }
\address{$^{3}$ Signal Processing Systems, Delft University of Technology, Delft, The Netherlands. {\it 
	E-mail address:} 
	{\sf{p.guerreroprado@tudelft.nl}}. }

\address{$^{4}$Departamento de Ingeniería Matemática and CI$^2$MA, Universidad de Concepción, Concepción, Chile. {\it 
	E-mail address:} 
	{\sf{fernandoroldan@udec.cl}}. }

\begin{document}
	\begin{abstract}
We introduce a primal-dual splitting algorithm with line search for monotone inclusions involving a
merely continuous operator, neither Lipschitzian nor monotone on the whole space, together with a
maximally monotone, a cocoercive, and a linearly composed maximally monotone operator. Its step-sizes
are selected via a new backtracking procedure which compares the variation of the continuous operator
with a joint primal-dual displacement. This condition is weaker than the purely primal tests used in
Tseng-type line searches, allowing for larger step-sizes. We show through examples that without the
dual part, the line search may fail to terminate or the iterates can converge to a point that is not a
solution. The line search is well defined under a local boundedness condition which is weaker than
those in the literature and cannot be removed. We prove weak convergence when either the step-sizes
are bounded away from zero or the continuous operator is uniformly continuous on weakly compact sets.
As an application, we consider convex problems involving the composition of a $K$-nondecreasing
function with a continuously differentiable $K$-convex mapping. Numerical experiments in computed
tomography with Poisson noise show that the proposed method reduces the computational time with
respect to existing methods with line search.   
		\par
		\bigskip
		
		
	\end{abstract}
	\maketitle
	\vspace*{-1cm}

    \section{Introduction}
	
    In this article, we aim to numerically solve the following problem.
\begin{pro}\label{pro:main}
	Let $\H$ and $\G$ be real Hilbert spaces, let $A: \H \to 2^\H$ and $B:\G \to 
	2^\G$ be maximally monotone operators, let $L\colon \H \to \G$ be a  bounded  linear operator with adjoint $L^*\colon \G \to \H$,
    let $X\subset \H$ be a nonempty closed convex set, let $C\colon \H \to \H$ be a $\beta$-cocoercive operator, for some $\beta \in 
	\RPP$, and let $D: \dom D \subset \H \to \H$ be a continuous operator on $\dom A\cup X \subset \dom D$. Moreover, assume that $A+D$ is maximally monotone. The problem is to
	\begin{equation} \label{eq:promain}
		\text{find}\quad (x,u) \in X \times \G \quad \text{such that}\quad
		\begin{cases}
			0 \in Ax+Cx+Dx+L^*u\\
			0\in B^{-1}u-Lx,
		\end{cases} 
	\end{equation}
	under the hypothesis that its set of solutions, denoted by 
	$\bm{Z}$,  is not 
	empty.
\end{pro}
 This monotone inclusion models several applied problems such as equilibrium problems \cite{Combettes2005equilibrium}, machine learning \cite{BotrelaxFBF2023,Nocedal2018,CombettesPesquet2021strategies,gidel2019variational,Mescheder2017}, optimization \cite{Combettes2018MP}, partial differential equations \cite{AubinHelene2009,Glowinsky1975,Showalter1997},  signal processing and imaging \cite{BotHendrich2014TV,Briceno2011ImRe,BurgerSawatzkySteidl2016,chambolle2016AN},  traffic theory \cite{Nets1,GafniBert84},
   variational inequalities \cite{bauschkebook2017},  among others. 

 In general, methods that rely on Lipschitzian operators require prior knowledge of the Lipschitz constant to initialize the iterations, and estimating this constant introduces an additional computational burden. Furthermore, the global step-size constraints dictated by these Lipschitz bounds can be far more conservative than the local bounds the algorithm could actually use during its execution. Moreover, relevant applications involve operators which are continuous but not Lipschitzian. This is the case in image denoising and reconstruction under Poisson noise \cite{Pustelnik2016}, where the data fidelity term is the Kullback--Leibler divergence, whose gradient is not Lipschitzian, and in convex optimization problems with nonlinear compositions \cite{MR4888251}, where the associated saddle operator is not Lipschitzian in general. These theoretical and practical limitations motivate the development of algorithms equipped with line search procedures.

In the particular case when $B=0$, $L=0$, $C=0$, $X=\H$, $A=\partial g$, and $D=\nabla f$, where $g$ is a proper lower semicontinuous convex function and $f$ is a differentiable convex function, Problem~\ref{pro:main} reduces to the minimization of $f+g$. When $\nabla f$ is Lipschitzian, this problem is solved by the forward-backward splitting \cite{Chenhg1997,combettes2005,passty1979JMAA}, whose step-sizes are bounded in terms of the Lipschitz constant. When this assumption does not hold, the convergence of the method is obtained in \cite{BelloCruz2016} and, in a variable metric setting, in \cite{Salzo2017}, by defining the step-sizes via line search procedures. Related line search procedures are proposed in \cite{BonettiniPratoRebegoldi2018} for problems satisfying the Kurdyka--\L ojasiewicz property, and in \cite{Bonettini2024} including an inertial term. Without line search, a step-size rule which adapts to the local behavior of the gradient is proposed in \cite{MalitskyMishchenko2024}, and the convergence of the forward-backward splitting with constant step-size is established in \cite{BricenoSilvaYang2024} when the gradient satisfies a locally Lipschitz-type condition.

More generally, when $C=0$, $B=0$, and $L=0$, Problem~\ref{pro:main} reduces to finding a zero of $A+D$ in $X$, which is studied in \cite{Tseng2000SIAM} for a general maximally monotone operator $A$ and an operator $D$ which is not necessarily a gradient. In this case the forward-backward splitting may fail to converge even if $D$ is Lipschitzian, and the method proposed there overcomes this difficulty by complementing the forward-backward step with an additional forward step, with a step-size bounded in terms of the Lipschitz constant of $D$. When $D$ is merely continuous, such a bound is not available and the step-size is found via a backtracking line search. 
More precisely, given the current iterate $x$ and $\theta\in\left]0,1\right[$, the step-size $\tau$ is reduced until
\begin{equation}\label{e:lsTseng}
	\tau\|Dx-Dz_x(\tau)\|\le\theta\|x-z_x(\tau)\|,
	\quad\text{where}\quad z_x(\tau)=J_{\tau A}(x-\tau Dx).
\end{equation}
Two features of this approach are relevant for our purposes. First, the step-size adapts to the local behavior of $D$, which allows for steps larger than those dictated by a global constant.
Second, $D$ need not be defined on the whole space and, although
\cite[Theorem~3.4]{Tseng2000SIAM} assumes $D$ to be maximally monotone,
its proof only uses the maximal monotonicity of $A$ and $A+D$. A variant of this method with an Armijo-type search requiring a single forward-backward computation per trial step is proposed in \cite{BelloCruzDiazMillan} in finite dimensions. Alternative strategies for this problem and related ones which avoid the knowledge of Lipschitz constants include the golden ratio algorithm \cite{Malitsky2020}, which adapts the step-size using local information in the context of variational inequalities, the line search variants \cite{MunozRoldan2026} of the forward-reflected-backward algorithm \cite{Malitsky2020SIAMJO}, and the line search for projective splitting methods in \cite{Eckstein2020}.

When $B=0$ and $L=0$ but $C\neq 0$, Problem~\ref{pro:main} reduces to finding a zero of $A+C+D$ in $X$. The method in \cite{Tseng2000SIAM} can be applied to this problem by considering the monotone continuous operator $C+D$, but $C$ is then evaluated twice per iteration and several additional times within the line search, and its cocoercivity is not exploited. This drawback is overcome by the forward-backward-half-forward (FBHF) splitting \cite{BricenoDavis2018}, in which $C$ is evaluated only once per iteration and the line search condition has the form \eqref{e:lsTseng} with $z_x(\tau)=J_{\tau A}(x-\tau(C+D)x)$, so that the test involves $D$ alone. As in \cite{Tseng2000SIAM}, its proof only uses the maximal monotonicity of $A+D$.

The approach in \cite{BricenoDavis2018} is extended in \cite{BricenoRoldan4op} to the problem
\begin{equation}\label{e:four}
	\text{find}\quad x\in X\quad\text{such that}\quad 0\in Ax+Cx+Dx+Ex,
\end{equation}
where $E\colon\H\to\H$ is a monotone and Lipschitzian operator. The algorithm proposed there evaluates $C$ once and $E$ twice per iteration, while the line search still involves only the continuous operator $D$. This setting covers Problem~\ref{pro:main} via a product space formulation. Indeed, \eqref{eq:promain} is the particular instance of \eqref{e:four} in $\H\times\G$ in which the maximally monotone operator is $(x,u)\mapsto Ax\times B^{-1}u$, the cocoercive operator is $(x,u)\mapsto(Cx,0)$, the continuous operator is $(x,u)\mapsto(Dx,0)$, and the Lipschitzian operator is the skew linear operator $(x,u)\mapsto(L^*u,-Lx)$. Hence, the algorithm in \cite{BricenoRoldan4op} provides a method with line search for Problem~\ref{pro:main}. The FBHF splitting \cite{BricenoDavis2018} can also be applied in $\H\times\G$ by considering the monotone continuous operator $(x,u)\mapsto(Dx+L^*u,-Lx)$, in which case $L$ and $L^*$ are additionally activated at each trial of the line search. However, in both approaches the linear operator is treated as a generic forward operator and the primal-dual structure of the problem is not exploited. In particular, the method in \cite{BricenoRoldan4op} requires two evaluations of $L$ and two of $L^*$ per iteration, and its step-sizes are subject to restrictive conditions.

On the other hand, when $D=0$, Problem~\ref{pro:main} is solved by the primal-dual algorithm in \cite{Vu13}, which is also proposed in \cite{Condat13} in the convex optimization setting and which extends the method in \cite{ChambollePock2011} to the case $C\neq 0$. The step-sizes of these methods are restricted by a condition involving $\|L\|$ and $\beta$. In the convex optimization setting, this restriction is removed in \cite{MalitskyPock2018}, where a line search for the algorithm in \cite{ChambollePock2011} is proposed in which each trial step requires updating only the dual variable. This line search is also extended to problems including an additional smooth function, and, in contrast with \cite{Vu13,Condat13}, the Lipschitz constant of its gradient is not needed. Related primal-dual line search rules are studied in \cite{ChangYangZhang2022} in combination with golden ratio techniques, with a preconditioned extension in \cite{ShanLiFeng2025}, in \cite{Nesterov2021} for linearly constrained minimization, and in \cite{Hamedani2021} for convex-concave saddle point problems with general coupling functions. A line search for a primal-dual method in the presence of a smooth function with locally Lipschitzian gradient is developed in \cite{AtenasDaoTam2025}, with applications to distributed optimization. Adaptive versions of the algorithm in \cite{Vu13,Condat13} which adjust the step-sizes to the local smoothness of the gradient are proposed in \cite{vladarean2021a,Latafat2025}. In the latter, no line search is needed for the gradient and a backtracking procedure is used only to avoid the computation of $\|L\|$.

In the case when $X=\H$ and $D$ is Lipschitzian, the primal-dual structure of Problem~\ref{pro:main} is fully exploited by the method proposed in \cite{Roldan20254op}. Given $(x_0,u_0)\in\H\times\G$, this method iterates as follows:
\begin{equation}\label{eq:algo1PD}
	(\forall n\in\N)\quad
	\begin{array}{l}
		\left\lfloor
		\begin{array}{l}
			z_{n+1} = J_{\tau A} \big(x_n-\tau (L^* u_n+Cx_{n}+Dx_n)\big) \\
			x_{n+1} = z_{n+1}-\tau(Dz_{n+1}-Dx_{n})\\
			u_{n+1} = J_{\sigma B^{-1}}\big(u_{n}+\sigma L(z_{n+1}+x_{n+1}-x_{n})\big).
		\end{array}
		\right.
	\end{array}
\end{equation}
When $D=0$, the recurrence in \eqref{eq:algo1PD} reduces to the algorithm in \cite{Vu13,Condat13}, and when $B=0$ and $L=0$ it reduces to FBHF \cite{BricenoDavis2018}. Hence, the algorithm in \cite{Roldan20254op} can be seen as a combination of these two methods. However, its step-sizes depend on the Lipschitz constant of $D$ and no line search is available when $D$ is merely continuous. Recently, a primal-dual forward-reflected splitting with line search is proposed in \cite{atenas2026primal} in finite-dimensional spaces for the case when $C=0$, $X=\H$, and $D$ is locally Lipschitzian, motivated by distributed min-max problems. In that method, the dual step uses the step-size of the previous iteration and the line search condition involves only the primal iterates, as in \eqref{e:lsTseng}.

To the best of our knowledge, none of the previous approaches provides a line search tailored to the primal-dual structure of Problem~\ref{pro:main} in its full generality, that is, in the simultaneous presence of the cocoercive operator $C$, a merely continuous operator $D$, and the linear composite term.
In this article, we fill this gap by proposing a primal-dual algorithm with line search for solving Problem~\ref{pro:main} in which the operator $D$ is only required to be continuous, so that no Lipschitz constant of $D$ needs to be known or estimated. The method is a variable step-size modification of \eqref{eq:algo1PD} in which the dual variable is updated first, the primal and dual variables share a single step-size, and the previous step-size enters the dual extrapolation. These modifications yield a Lyapunov function that decreases under suitable conditions on the step-size, which is defined via a new primal-dual line search condition of the form
\begin{equation*}
	\tau^2\|Dx-Dz_x(\tau)\|^2\le\theta^2\big(\|x-z_x(\tau)\|^2+\|R(\tau)\|^2\big),
\end{equation*}
where $x$ is the current primal iterate, $z_x(\tau)$ is the outcome of the primal forward-backward step, and $R(\tau)$ is a residual involving the dual iterates.
In other words, the line search compares the variation of $D$ along the trial step with the joint primal-dual displacement. We justify the inclusion of this primal-dual displacement measure by presenting an example in which a purely primal line search, as in \eqref{e:lsTseng}, fails to terminate (see Example~\ref{ex:primalLS}). Furthermore, since this displacement includes a nonnegative dual term in addition to the primal residual in \eqref{e:lsTseng}, the condition is weaker than the purely primal tests used in \cite{Tseng2000SIAM,BricenoDavis2018,atenas2026primal}, which allows for larger step-sizes, as shown by our numerical tests. In addition, when $D$ is not Lipschitzian, this dual term is essential for convergence, as illustrated in Example~\ref{ex:primalLSconvL} where the primal line search terminates but the iterates converge to a point which is
not a solution. In exchange, each trial of the line search requires one evaluation of the resolvent of $B^{-1}$ and one of $L^*$. Moreover, the operator $C$ is evaluated only once per iteration and it is not activated within the line search. We show that the line search procedure terminates whenever $D$ is locally bounded on
$\overline{\dom}A$ relative to $\dom A$. This condition is implied by the assumptions guaranteeing the termination of the line searches in \cite{Tseng2000SIAM,BricenoDavis2018,BricenoRoldan4op}
(see Remark~\ref{rem:LS}) and, as shown in Example~\ref{ex:notbound}, it cannot be removed, even if $D$ is maximally monotone. Furthermore, we derive weak convergence of the primal and dual iterates when either the accepted step-sizes remain bounded away from zero or $D$ is uniformly continuous on weakly compact sets. When $D$ is Lipschitzian and the step-size is constant, the method and its step-size condition reduce to those of \cite{Roldan20254op}. Moreover, in contrast with \cite{Roldan20254op}, the set $X$ is not required to be the whole space. Unlike the line search methods in \cite{BricenoDavis2018,BricenoRoldan4op}, which solve Problem~\ref{pro:main} in the primal-dual product space and treat the skew operator $(x,u)\mapsto(L^*u,-Lx)$ as an additional forward operator, our scheme handles $L$ as in \cite{Vu13,Condat13}.

Although \cite{Tseng2000SIAM,BricenoDavis2018,BricenoRoldan4op} assume $D$ to be maximally monotone, their proofs only use the maximal monotonicity of $A+D$. Accordingly, our convergence result does not require $D$ to be maximally monotone or monotone on the whole space. The only requirement is the maximal monotonicity of the sum $A+D$. This feature is essential in convex optimization problems involving nonlinear compositions \cite{MR4888251}, such as
\begin{equation}\label{e:nonlincomp}
	\minimize{x\in\H}{f(x)+h(x)+g(Lx)+\psi(d(x))},
\end{equation}
where $f$, $g$, and $\psi$ are proper lower semicontinuous convex functions, $h$ is a differentiable convex function with a Lipschitzian gradient, $L$ is a bounded linear operator, $d$ is a continuously differentiable $K$-convex mapping for some closed convex cone $K$, and $\psi$ is $K$-nondecreasing. This setting includes, for instance, the pointwise maximum of finitely many differentiable convex functions $d_1,\ldots,d_m$, which is obtained when $K$ is the nonnegative orthant of $\R^m$, $d=(d_1,\ldots,d_m)$, and $\psi\colon y\mapsto\max_{1\le i\le m}y_i$. Such problems arise in worst-case robust approximation \cite[Section~6.4.2]{BoydVandenberghe2004}, in distributionally robust and federated learning \cite{Sagawa2020,Mohri2019}, and in discrete minimax theory \cite{DemyanovMalozemov1974}. The primal-dual formulation of \eqref{e:nonlincomp} is a particular instance of Problem~\ref{pro:main} in which the continuous operator is $(x,w)\mapsto\big([\nabla d(x)]^*w,-d(x)\big)$, where $w$ is the dual variable associated with the nonlinear composition. This operator is, in general, neither Lipschitzian nor monotone on the whole space, but it is monotone on the domain of the corresponding set-valued operator, and their sum is maximally monotone. Hence, the methods relying on the Lipschitz continuity or on the maximal monotonicity of each operator are not applicable, while our algorithm and its line search exploit the full structure of \eqref{e:nonlincomp}, as detailed in Section~\ref{sec:aplication}.

To illustrate the applicability of our algorithm, we perform numerical experiments on total variation minimization under a Kullback--Leibler discrepancy constraint, arising in computed tomography reconstruction with Poisson noise.
We compare our method with the natural competitors from \cite{BricenoDavis2018,BricenoRoldan4op}. The results show that the proposed algorithm reduces the average computational time and the number of iterations with respect to these methods.

The remainder of the paper is organized as follows. In Section~\ref{sec:notation}, we introduce our notation and preliminary results. In Section~\ref{sec:main}, we present our algorithm and its convergence analysis. Section~\ref{sec:aplication} is devoted to developing applications in optimization problems. In Section~\ref{sec:num}, we present our numerical experiments. Finally, Section~\ref{sec:conclu} provides our concluding remarks. The details of the examples in Section~\ref{sec:main} are provided in the Appendix.

\section{Notation and Preliminaries}\label{sec:notation}
Throughout this paper, $\H$ and $\G$ denote real Hilbert spaces equipped with the inner product $\scal{\cdot}{\cdot}$ and associated norm $\|\cdot \|=\sqrt{\scal{\cdot}{\cdot}}$. For a given $x \in \H$ and a radius $\delta \in \RPP$, the open ball of radius $\delta$ centered at $x$ is written as $\mathcal{B}(x,\delta)$. We use $\to$ to indicate strong convergence, whereas $\weak$ stands for weak convergence. The identity operator on a given space is represented by $\id$. We frequently rely on two fundamental inequalities:
\begin{eqnarray}\label{eq:CS}
	&\hspace*{-1cm}(\textnormal{Cauchy–Schwarz inequality})  &(\forall (x,u) \in \H^2) \quad |\scal{x}{u}| \leq \|x\|\|u\|,\end{eqnarray}
\begin{eqnarray}
	\label{eq:YI}&(\textnormal{Young's inequality})   &(\forall (a,b)\in \R^2)(\forall \epsilon >0) \quad 2ab \leq \epsilon a^2+\frac{b^2}{\epsilon}.
\end{eqnarray}

Consider a nonempty subset $D\subset \H$, an operator $T\colon D \rightarrow \H$, and a constant $\beta \in \RPP$. The operator $T$ is $\beta$-cocoercive if
\begin{equation*} 
	(\forall (x,y) \in D^2) \quad \scal{x-y}{Tx-Ty} \geq \beta \|Tx - Ty\|^2.
\end{equation*}
Furthermore, $T$ is said to be $\beta$-Lipschitz continuous on $D$ if
\begin{equation*} 
	(\forall (x,y) \in D^2) \quad \|Tx-Ty\| \leq \beta\|x - y \|.
\end{equation*}
It follows from the Cauchy–Schwarz inequality that any $\beta$-cocoercive operator is $(1/\beta)$-Lipschitz continuous. The operator $T$ is continuous on $D$ if for every $u \in D$ and any sequence $(u_n)_{n \in \N} \subset D$ such that $\|u_n - u\| \to 0$, we have $\|Tu_n - Tu\| \to 0$. On the other hand, $T$ is uniformly continuous on $D$ if for every pair of sequences $(u_n)_{n \in \N}$ and $(v_n)_{n \in \N}$ in $D$ satisfying $\|u_n - v_n\| \to 0$, it holds that $\|Tu_n - Tv_n\| \to 0$.

The power set of $\H$ is denoted by $2^\H$. For a set-valued mapping $A\colon\H \rightarrow 2^{\H}$, its domain, 
graph, and the set of zeros of $A$ are, respectively,
$\dom A=\menge{x\in\H}{Ax\neq\varnothing}$, 
$\zer A=\menge{x\in\H}{0\in A x}$, and $\gra A=\menge{(x,u)\in\H\times\H}{u\in A x}$.
The inverse operator $A^{-1}\colon\H \rightarrow 2^{\H}$ is given by $u \mapsto \menge{x \in \H}{u \in Ax}$.
Recall that $A$ is monotone if 
\begin{equation*}
	\left(\forall \big((x,u),(y,v)\big) \in (\gra A)^2\right) \quad \scal{x-y}{u-v} \geq 0.
\end{equation*}
Moreover, $A$ is maximally monotone if there exists no monotone operator whose graph strictly contains $\gra A$. If $A$ is a maximally monotone operator, $\overline{\dom} A$ is convex. The resolvent associated with $A$ is the single-valued mapping defined as $J_A\colon\H \rightarrow \H \colon x \mapsto (\id+A)^{-1}x$. 

 Let $f\colon\H\to\RX$. The domain 
of $f$ is $\dom f=\menge{x\in\H}{f(x)<+\infty}$, and $f$ is proper if 
$\dom f\neq\varnothing$. We denote by $\Gamma_0(\H)$ the class of all proper, 
lower semicontinuous, convex functions from $\H$ to $\RX$. Let 
$f\in\Gamma_0(\H)$. The Fenchel conjugate of $f$ is the function 
$f^*\in\Gamma_0(\H)$ defined by 
$f^*\colon u\mapsto \sup_{x\in\H} \left(\scal{x}{u}-f(x)\right)$.
The subdifferential of $f$ is the operator $\partial f\colon x\mapsto \menge{u\in\H}{(\forall y\in\H)\:\: f(x)+\scal{y-x}{u}\le f(y)}$,
it is maximally monotone, it satisfies 
$(\partial f)^{-1}=\partial f^*$, and its set of zeros $\zer\partial f$ 
coincides with the set of global minimizers of $f$, denoted by 
$\argmin_{x\in \H}f$. If $f$ is Fr\'echet differentiable at $x\in\H$, then 
$\partial f(x)=\{\nabla f(x)\}$. The proximity operator of $f \in \Gamma_0(\H)$ is the single-valued mapping
\begin{equation*}
	\prox_{f}\colon x\mapsto\argmin_{y\in\H}\left(f(y)+\frac{1}{2}\|x-y\|^2\right),
\end{equation*}
which precisely corresponds to the resolvent $J_{\partial f}$. Given a nonempty closed convex set $S\subset\H$, $P_S$, $N_S$, and $\iota_S$ denote the projection
onto $S$, the normal cone to $S$, and the indicator function of $S$, respectively, and $\sri S$ denotes the strong relative interior of $S$ \cite[Definition~6.9]{bauschkebook2017}. For a comprehensive treatment of monotone operator theory and convex analysis, we direct the reader to \cite{bauschkebook2017}. We conclude this section by introducing the following lemma, a direct consequence of \cite[Proposition 3.4]{MR5098995} and \cite[Lemma~3.2]{BricenoRoldan4op}.
\begin{lem}\label{lemma}
		Let $A\colon \H \to 2^\H$ be a maximally monotone operator, 
    let $(z,y) \in \H^2$, and define 
	\begin{equation*}
		(\forall \tau > 0) \quad \quad	x_{z,y}(\tau) = 
		J_{\tau A} 
		( z - 
		\tau y). 
	\end{equation*}
	Then, the following assertions hold:
	\begin{enumerate}
    \item\label{lemma0} $x_{z,y}(\tau) \to P_{\overline{\dom} A}z$ as $\tau \downarrow  0$.
		\item\label{lemmai} The function 
		$\tau \mapsto \frac{1}{\tau} \|z- x_{z,y}(\tau)\|$
		is nonincreasing.
	\end{enumerate}
\end{lem}

In Problem~\ref{pro:main}, we assume that $A+D$ is maximally monotone, which does not
require $D$ to be monotone on the whole space. In the next result, we provide two instances in which the aforementioned assumption holds.
\begin{lem}\label{l:A+D}
Let $A\colon\H\to2^\H$ be maximally monotone and let $D\colon\dom D\subset\H\to\H$ with $\dom A\subset\dom D$. Then $A+D$ is maximally monotone if one of the following holds.
\begin{enumerate}[label=\textup{(\alph*)}]
\item\label{it:convex}    $D$ is monotone and
continuous on a nonempty
closed convex set  $S\subset\H$ such that $\dom A\subset S\subset\dom D$ and
$\dom A\cap\inte S\neq\varnothing$.

\item\label{it:saddle} $\H=\H_1\times\H_2$ and $A+D$ is the saddle subdifferential operator $(x,w)\mapsto\partial(\mathcal L(\cdot,w))(x)\times\partial(-\mathcal L(x,\cdot))(w)$ of a closed proper convex--concave function $\mathcal L$.
\end{enumerate}    
\end{lem}
\begin{proof}
\ref{it:convex} By \cite[Theorem~3]{MR282272}, $N_S+D\vert_S$ is maximally monotone with domain $S$.
Since $\dom A\subset S$, we have $0\in N_Sx$ for every $x\in\dom A$, so
$\gra A\subset\gra(A+N_S)$, and the maximality of $A$ yields $A+N_S=A$.
Hence $A+D=A+(N_S+D\vert_S)$, which is maximally monotone by
\cite[Corollary~25.5]{bauschkebook2017}, because
$\dom A\cap\inte\dom(N_S+D\vert_S)=\dom A\cap\inte S\neq\varnothing$.

\ref{it:saddle} This is Rockafellar's theorem on saddle functions; see, e.g., \cite{Roc70b} and \cite[Theorem~4.3]{Combettes2018MP}.
 \end{proof}

\section{Main Results}\label{sec:main}
In this article, we propose the following algorithm for solving 
Problem~\ref{pro:main}.
\begin{algo}\label{algo:algo1}
In the context of Problem~\ref{pro:main}, let $(x_{-1},x_0,z_0,u_0)\in 
\H\times X \times \H \times \G$, let $\tau\in \RPP$, let $\tau_{-1}\in\ ]0,\tau]$, let $(\tau_n)_{n \in \N}$ be a sequence in $]0,\tau]$, and consider the recurrence
\begin{equation}\label{eq:algo1}
(\forall n\in\N)\quad 
\begin{array}{l}
\left\lfloor
\begin{array}{l}
v_n = z_n-x_{n-1}\\
u_{n+1} = J_{\tau_n B^{-1}} \big(u_n+\tau_nLx_n+\tau_{n-1}Lv_n\big)\\
z_{n+1} = J_{\tau_n A} \big(x_n-\tau_n (L^* u_{n+1}+Cx_n+Dx_n)\big) \\			
x_{n+1} =P_X \big(z_{n+1}-\tau_n(Dz_{n+1}-Dx_n)\big).
\end{array}
\right.
\end{array}
\end{equation}
\end{algo}   
Unlike \eqref{eq:algo1PD}, Algorithm~\ref{algo:algo1} updates the dual variable first and includes the
previous step-size $\tau_{n-1}$ in the dual update, which yields the Lyapunov-type estimate for
iteration-dependent step-sizes in Proposition~\ref{prop:des}. In order to define the step-sizes $(\tau_n)_{n \in \N}$ in Algorithm~\ref{algo:algo1}, we generalize the procedure that defines the step-sizes for the three-operator splitting method in \cite{BricenoDavis2018}. Our test measures the increment in evaluating $D$ with a primal-dual residual, in contrast to traditional line search conditions that only use primal displacements.

\begin{line search}\label{line search1}
		In the context of Problem~\ref{pro:main}, let $\x=(x,z,v,u) \in  X  \times \H^2 \times \G$, $ \tau \in \RPP$, $(\theta,\,\gamma) \in \left]0,1\right[^2$, and $\tilde{\tau} \in \RPP$. The line search returns $\alpha_k = \tau \gamma^k$ with $k \in \N$ being the smallest natural number such that 
		\begin{equation}\label{eq:line search1}
			\alpha_k^2\|D\fz(\alpha_k,\x)-Dx\|^2\leq \theta^2 (\|\fz(\alpha_k,\x)-x\|^2 + \|u-\fu(\alpha_k,\x)+\tilde{\tau}Lv \|^2),
		\end{equation}
		where 
	\begin{equation}\label{eq:defqppk}
		\begin{cases}
			\fu(\alpha,\x) = J_{\alpha B^{-1}} \big(u+\alpha Lx+\tilde{\tau}Lv\big)\\
            \fz(\alpha,\x) = J_{\alpha A} \big(x-\alpha (L^* \fu(\alpha,\x)+Cx+Dx)\big).	
	\end{cases}	
	\end{equation}
\end{line search}	

The following example shows that, indeed, the term $\theta^2 \|u-\fu(\alpha_k,\boldsymbol{x})+\tilde{\tau} L v\|^2$ cannot be dropped from \eqref{eq:line search1}  while still guaranteeing that the line search terminates. 
\begin{exam}\label{ex:primalLS}
In the context of Problem~\ref{pro:main}, set $\H=\G=\R$, $X=\R$, $A=0$, $C=0$, $L=\id$, $B=\id$, and $	D\colon\R\to\R\colon x\mapsto x^{1/3}+1$. All the assumptions of Problem~\ref{pro:main} are verified (see the Appendix) and $\bm{Z}=\{(\bar x,\bar x)\}$ with $\bar x = t^3$, where $t$ is the unique real solution to the equation $t^3 + t + 1 = 0$. Consider $\x=(x,z,v,u)=(0,0,0,-1)$ and $\tilde\tau\in\RPP$. In this setting, if the dual residual is removed from \eqref{eq:line search1}, the test becomes
$\alpha^2|D\fz(\alpha,\x)-Dx|^2\le\theta^2|\fz(\alpha,\x)-x|^2$, which is equivalent to  $(1+\alpha)^{4/3}\le\theta^2\alpha^{2/3}$. Note that, since $\theta \in ]0,1[$, this inequality fails for
every $\alpha\in\RPP$. Hence, the
line search with the purely primal test never stops. On the other hand, the test in
\eqref{eq:line search1} reduces to $\big(\alpha(1+\alpha)\big)^{4/3}\le\theta^2(1+\alpha^2)$,
which holds whenever $\alpha(1+\alpha)\le\theta^{3/2}$. Since $\alpha_k=\tau\gamma^k\to0$, Line
Search~\ref{line search1} terminates. 
\end{exam}

The following proposition asserts that Line Search \ref{line search1} is well defined.
\begin{prop}\label{prop:LS}
In the context of Line Search~\ref{line search1}, denote by 
$\x=(x,z,v,u)\in X  \times \H^2 \times \G$ and suppose that 
one of the following holds.
\begin{enumerate}
    \item\label{prop:LS1} $x\in {\overline{\dom}A}$.
    \item\label{prop:LS2} $x\notin {\overline{\dom}A}$ and there exists $\rho>0$ such that $D\big(\mathcal{B}(P_{\overline{\dom}A}x,\rho)\cap\dom A\big)$ is bounded.
\end{enumerate}
Then, Line Search~\ref{line search1} terminates. 
\end{prop}	
\begin{proof} Suppose that Line Search~\ref{line search1} does not terminate. Then $\alpha_k \to 0$ and \begin{equation}\label{eq:propLS1}
	(\forall k \in \N)	\quad \alpha_k^2\|Dp_k-Dx\|^2> \theta^2 (\|p_k-x\|^2 + 
    \| u-q_k+\tilde{\tau}Lv\|^2),
	\end{equation} where $p_k:=\fz(\alpha_k,\boldsymbol{x})$ and $q_k:=\fu(\alpha_k,\boldsymbol{x})$.
    Set $w=P_{\overline{\dom}B^{-1}}(u+\tilde{\tau}Lv)$ and $y=P_{\overline{\dom}A}(x)$. It follows from Lemma~\ref{lemma}\eqref{lemma0} and \eqref{eq:defqppk} that $q_k \to w$ and, in turn, that $p_k\to y$. Indeed, for each $k \in \N$, define $\widetilde{p}_{k} = J_{\alpha_k A} \big(x-\alpha_k (L^* w+Cx+Dx)\big)$. Then, $\widetilde{p}_{k} \to y$  in view of Lemma~\ref{lemma}\eqref{lemma0} and,
by the nonexpansiveness of $J_{\alpha_k A}$,
$ \| \widetilde{p}_{k} - p_k  \| \leq \alpha_k \| L^* (w-q_k) \| \to 0$.

\ref{prop:LS1}: In this case $y=x$ and 
    $p_k \to x$. Since $\{p_k\}_{k \in \N}\cup \{x\} \subset \dom A\cup X$, by the continuity of $D$ in $\dom A\cup X$ we have that $Dp_k \to Dx$. Therefore, \eqref{eq:propLS1} yields
 \begin{align}\label{eq:propLS2}
 \frac{x-p_k}{\alpha_k} \to 0 \quad \textnormal{ and } \quad  \frac{u-q_k+\tilde{\tau}Lv}{\alpha_k} \to 0.
 \end{align}
Moreover, it follows from \eqref{eq:defqppk} that
\begin{align}\label{eq:propLS3}
			\dfrac{x-p_k}{\alpha_k } -(L^*q_k+Cx+Dx) \in A p_k \quad \textnormal{ and }\quad 
			\dfrac{u-q_k+\tilde{\tau}Lv}{\alpha_k } + L x \in B^{-1} q_k.
\end{align}
Note that \eqref{eq:propLS2} yields $(q_k,p_k) \to  (u+\tilde{\tau}Lv,x)$ and since $A$ and $B$ are maximally monotone operators, by \eqref{eq:propLS2}, \eqref{eq:propLS3}, and the closure of the graph of maximally monotone operators \cite[Proposition~20.38]{bauschkebook2017}, we deduce that
	\begin{align} \label{eq:solLS}
		0 \in Ax+Cx+Dx+L^*(u+\tilde{\tau}Lv) \quad \textnormal{ and } \quad
		0 \in B^{-1}(u+\tilde{\tau}Lv)-Lx.
\end{align}
Then $(x,u+\tilde{\tau}Lv) \in \bm{Z}$. In this case, \eqref{eq:solLS} is equivalent to
	\begin{align*}
    \begin{cases}
	x = J_{\tau A} ( x-\tau(Cx+Dx+L^*(u+\tilde{\tau}Lv)))\\
	u+\tilde{\tau}Lv =J_{\tau B^{-1}} (u+\tilde{\tau}Lv+	\tau Lx) 
	\end{cases} \Leftrightarrow 
	\begin{cases}
	p_0 =	x = J_{\tau A} ( x-\tau(Cx+Dx+L^*q_0))\\
	q_0 =	u+\tilde{\tau}Lv =J_{\tau B^{-1}} (u+\tilde{\tau}Lv+	\tau Lx). 
	\end{cases}
	\end{align*} 
	Therefore,  \eqref{eq:line search1} holds for $k=0$, a contradiction with \eqref{eq:propLS1}.

\ref{prop:LS2}: Let $\delta\in\left]0,\min\{\rho,\|y-x\|/2\}\right[$. Since $p_k \to y$, there exists
$k_0 \in \mathbb{N}$ such that, for every $k \geq k_0$, $\|p_k-y\|<\delta$ and, hence,
$\|p_k-x\|\ge\|y-x\|-\|p_k-y\|>\delta$ and $Dp_k\in D(\mathcal{B}(y, \delta) \cap \dom A)\subset
D(\mathcal{B}(y, \rho) \cap \dom A)$.
Therefore, there exists $M>0$ such that $\sup_{k\ge k_0}\|Dp_k-Dx\|\le M$.
Moreover, since $\alpha_k \to 0$, there exists 
$k_1\ge k_0$ such that, for every $k\ge k_1$, $\alpha_k^2\le \theta^2\delta^2/M^2$ leading to
\begin{equation*}
	\alpha_k^2\|Dp_k-Dx\|^2\leq \theta^2\|p_k-x\|^2\leq \theta^2(\|p_k-x\|^2+\|  u-q_k+\tilde{\tau}Lv\|^2),
\end{equation*} yielding a contradiction with \eqref{eq:propLS1}, and the proof is complete.
 \end{proof}

\begin{rem}\label{rem:LS}
\begin{enumerate}
\item If $X\subset\overline{\dom}A$, assertion \ref{prop:LS1} in Proposition~\ref{prop:LS} holds for every
$x\in X$ and Line Search~\ref{line search1} always terminates.
\item Assumption \ref{prop:LS2}  in Proposition~\ref{prop:LS} holds for every $x\in X\setminus\overline{\dom}A$
provided that $D\vert_{\dom A}$ is locally bounded at every point of
$\overline{\dom}A$, i.e.,
\begin{equation}\label{eq:locbound}
(\forall y\in\overline{\dom}A)(\exists\rho>0)\quad
D\big(\mathcal{B}(y,\rho)\cap\dom A\big)\ \text{is bounded}.
\end{equation}
Since $D$ is continuous on $\dom A\cup X$, \eqref{eq:locbound} holds
automatically at every $y\in\overline{\dom}A\cap(\dom A\cup X)$. Hence
\eqref{eq:locbound} is a genuine restriction only at points of
$\overline{\dom}A\setminus(\dom A\cup X)$ and, in particular, it
holds in each of the following cases.
\begin{enumerate}
\item $\overline{\dom}A\subset\dom A\cup X$, e.g., $\dom A$ is closed.
\item $\overline{\dom}A\subset\dom D$ and $D\vert_{\overline{\dom}A}$ is continuous.
This is the case, for instance, if $D$ is continuous on an open set
containing $\overline{\dom}A$.
\item $\overline{\dom}A\subset\inte\dom D$ and $D$ is maximally monotone
\cite[Theorem~1]{Rockafellar1969}.
\item  $D$ is bounded on bounded subsets of $\dom A$, e.g., $D$ is
Lipschitz continuous on $\dom A$.
\end{enumerate}
In Example~\ref{ex:notbound} we provide an instance in which \eqref{eq:locbound} fails, $(Dp_k)_{k\in\N}$ is unbounded, and Line Search~\ref{line search1} does not terminate.
\item In our notation, the termination of the line search in
\cite[Theorem~3.4(a)]{Tseng2000SIAM} is guaranteed when $D$ is Lipschitz
continuous on $X\cup\dom A$, or $D$ is continuous on $\dom A$ and
$X\subset\dom A$. In \cite[Theorem~2.3]{BricenoDavis2018} it is guaranteed
under the assumption $X\subset\dom A$, and in
\cite[Theorem~3.4]{BricenoRoldan4op} when $X\subset\dom A$, or
$\overline{\dom}A\subset\dom D$ and $D$ is locally bounded in $\dom D$, or
$\overline{\dom}A\subset\inte\dom D$. Each of these conditions implies
\ref{prop:LS1} or \eqref{eq:locbound}. Moreover, in \cite{BricenoDavis2018,BricenoRoldan4op,Tseng2000SIAM} the authors assume
that $D$ is maximally monotone, which is not required in
Proposition~\ref{prop:LS}.

    \item
    It follows from the proof of Proposition~\ref{prop:LS} that if $(x,u+\tilde{\tau}Lv)\in\bm{Z}$, Line Search~\ref{line search1} terminates at $k=0$, and $(\fz(\tau,\boldsymbol{x}),\fu(\tau,\boldsymbol{x}))=(x,u+\tilde{\tau}Lv)\in\bm{Z}$.
\end{enumerate}
\end{rem}
\begin{exam}\label{ex:notbound}
In the context of Problem~\ref{pro:main}, set $\H=\R^2$, $\G=\R$, $L=0$, $B=0$,
$X=\menge{(t,r)\in\R^2}{t+2r\ge1}$, $C\colon\R^2\to\R^2\colon(t,r)\mapsto(1,0)$, $A=\partial f$, and $D = \partial d$, where
\begin{equation*}
	f\colon(t,r)\mapsto
	\begin{cases}
		-\ln t, & \text{if } t>0 \text{ and } r=0,\\
		\pinf, & \text{otherwise}, 
	\end{cases}\quad \textnormal{and} \quad d\colon(t,r)\mapsto
	\begin{cases}
		-2\ln (t+2r), & \text{if } t+2r>0,\\
		\pinf, & \text{otherwise}.
	\end{cases}
\end{equation*} We have that
$\bm{Z}=\{((3,0),0)\}$. 
Let $\x=(x,z,v,u)$ with $x=(-1,1)\in X$, $u=0$, and $(z,v)\in\R^2\times\R^2$, let
$\tilde\tau\in\RPP$, let $\tau\in\left]0,1/2\right]$, and let $\gamma \in ]0,1[$. In this setting, the test
\eqref{eq:line search1} fails for every $k\in\N$, and Line Search~\ref{line search1} does not terminate. Moreover, neither \ref{prop:LS1} nor \ref{prop:LS2} in Proposition~\ref{prop:LS} holds and, indeed,
$\|D\fz(\alpha_k,\x)\|\to\pinf$. Note also that $D$ is maximally monotone, and thus this
property alone does not guarantee the termination of the line search. For additional details see the Appendix.
\end{exam}

\begin{prop}\label{prop:des}
	In the context of Problem~\ref{pro:main}, let $(x_{-1},x_0,z_0,u_0)\in 
	\H\times X \times \H \times \G$, and let 
	$(z_n)_{n \in 
		\N}$, $(x_n)_{n \in \N}$, and $(u_n)_{n \in \N}$ be the 
	sequences generated by Algorithm~\ref{algo:algo1}. Let $\theta\in\left]0,1\right[$ and suppose that
    \begin{equation}\label{eq:propdestau}
    \tau \in \left]0,\chi(\theta)\right[, \quad \textnormal{ where } \quad \chi(\theta)=\frac{4\beta(1-\theta^2)}{1+\sqrt{1+16\beta^2\|L\|^2(1-\theta^2)}}.
        \end{equation}
    For every $n \in \N$ and $(x,u) \in \bm{Z}$, define
	\begin{equation}
\label{def:gamman}
        \Gamma_{n}(x,u):=\|x_{n}-x\|^2 + \|u_{n}-u\|^2 + \left(1-\frac{\tau}{2\beta}-\theta^2\right)\|v_n\|^2+2\tau_{n-1}\scal{Lv_n}{u_{n}-u}.
	\end{equation}
Suppose that 
		\begin{equation}\label{eq:desls1}
(\forall n \in \N)\quad  \tau_n^2\|Dz_{n+1}-Dx_n\|^2\leq \theta^2 (\|v_{n+1}\|^2 + \| u_n-u_{n+1}+\tau_{n-1}Lv_n \|^2).
		\end{equation}
	Then, for any  $(x,u) \in \bm{Z}$, the following assertions hold:
	\begin{enumerate}
		\item 		\label{prop:des1} 
 For every $n \in \N$, 
		\begin{align}\label{eq:propdes1} 
        		\Gamma_{n+1}(x,u) \leq \Gamma_n (x,u) -
        \left(1-\frac{\tau}{2\beta}-\theta^2-\tau^2\|L\|^2\right)  \|v_n\|^2& 
    -(1-\theta^2)\|u_n-u_{n+1}+\tau_{n-1}Lv_n\|^2.
		\end{align}		
        \item\label{prop:des15}  $(\Gamma_{n}(x,u))_
        {n \in \N}$ is nonincreasing.
		\item \label{prop:des2}	For every $n \in \N$,
		\begin{equation*}
				\Gamma_{n}(x,u)
			\geq \|x_{n}-x\|^2 + \left(1-	\frac{\tau^2\|L\|^2}{(1-\tau/(2\beta)-\theta^2)} \right)\|u_{n}-u\|^2\geq 0.
		\end{equation*}
	\end{enumerate}
\end{prop}	
	\begin{proof}
        \eqref{prop:des1} Let $(x,u) \in \bm{Z}$. It follows from \eqref{eq:algo1} that
\begin{equation*} 
	\frac{u_n-u_{n+1}}{\tau_n}+Lx_n+ \frac{\tau_{n-1}}{\tau_n}Lv_n \in B^{-1}u_{n+1}.
\end{equation*}
By invoking the monotonicity of $B^{-1}$, we obtain
\begin{equation*} 
0\leq 	\scal{\frac{u_n-u_{n+1}}{\tau_n}+Lx_n+ \frac{\tau_{n-1}}{\tau_n}Lv_n-Lx}{u_{n+1}-u}.
\end{equation*}	
Using the polarization identity,
\begin{align*}
	\frac{1}{\tau_n}\scal{u_n-u_{n+1}}
	{u_{n+1}-u}
	=\frac{1}{2\tau_n}(\|u_n-u\|^2-\|u_n-u_{n+1}\|^2-\|u_{n+1}-u\|^2),
\end{align*}
the previous inequality becomes
\begin{equation}\label{eq:proof1}
\|u_{n+1}-u\|^2\leq	\|u_n-u\|^2-\|u_n-u_{n+1}\|^2+2\tau_n\scal{Lx_n-Lx}{u_{n+1}-u}+2\tau_{n-1}\scal{Lv_n}{u_{n+1}-u}.
\end{equation}Analogously, the third line in \eqref{eq:algo1} gives
\begin{equation*}
\frac{x_n-z_{n+1}}{\tau_n}-(L^*u_{n+1}+Cx_n+Dx_n)+Dz_{n+1}\in (A+D)z_{n+1}.
\end{equation*}
Since $-(L^*u+Cx)\in(A+D)x$ by \eqref{eq:promain}, the monotonicity of $A+D$ yields
\begin{equation*}
0\le\scal{\frac{x_n-z_{n+1}}{\tau_n}-L^*(u_{n+1}-u)-(Cx_n-Cx)-(Dx_n-Dz_{n+1})}{z_{n+1}-x}.
\end{equation*}
Again, by the polarization identity,
we conclude that 
\begin{align}
\label{eq:proof2} \scalebox{0.98}{$
\|z_{n+1}-x\|^2\leq \|x_n-x\|^2-\|v_{n+1}\|^2-2\tau_n\scal{L^*(u_{n+1}-u)+Cx_n-Cx+Dx_n-Dz_{n+1}}{z_{n+1}-x}$}.    
\end{align}
Summing \eqref{eq:proof1} and \eqref{eq:proof2} provides the estimate
\begin{align}\label{eq:eq:proof3}
\|z_{n+1}-x\|^2 + \|u_{n+1}-u\|^2 &+  \|v_{n+1}\|^2+2\tau_n\scal{Lv_{n+1}}{u_{n+1}-u} \nonumber\\
& \leq  \|x_{n}-x\|^2 + \|u_{n}-u\|^2 - \|u_{n+1}-u_n\|^2 +2\tau_{n-1}\scal{Lv_n}{u_{n+1}-u}\nonumber\\
&\qquad\quad -2\tau_{n}\scal{Cx_n-Cx}{z_{n+1}-x} -2\tau_{n}\scal{Dx_n-Dz_{n+1}}{z_{n+1}-x}.
\end{align}
On the other hand, since $x \in X$, the definition of $x_{n+1}$ in \eqref{eq:algo1} implies
	\begin{align}
	 \|x_{n+1}-x\|^2 & =  \|P_X(z_{n+1}-\tau_n(Dz_{n+1}-Dx_n))-x\|^2 \nonumber\\
	 &\leq  
	 \|z_{n+1}-\tau_n(Dz_{n+1}-Dx_n)-x\|^2 \nonumber\\
	&= \|z_{n+1}-x\|^2-2\tau_n\scal{Dz_{n+1}-Dx_n}{z_{n+1}-x}+\tau_n^2\|Dz_{n+1}-Dx_n\|^2.\label{eq:eq:proof4}
	\end{align}
Furthermore, from the cocoercivity of $C$, \eqref{eq:CS}, and \eqref{eq:YI} with $\epsilon = 2\beta$,  we have
\begin{align}\label{eq:eq:proof5}
			2\tau_n\scal{Cx_n-Cx}{x-z_{n+1}} &=  
			2\tau_n\scal{Cx_n-Cx}{x-x_{n}}+2\tau_n\scal{Cx_n-Cx}{x_{n}-z_{n+1}}\nonumber\\
			&\leq 
			-2\tau_n\beta\|Cx_n-Cx\|^2+2\tau_n\scal{Cx_n-Cx}{x_{n}-z_{n+1}}\nonumber\\
            &\le \frac{\tau_n}{2\beta}
			\|v_{n+1}\|^2\nonumber\\
            &\leq \frac{\tau}{2\beta}
			\|v_{n+1}\|^2.
		\end{align}
Combining \eqref{eq:eq:proof3}, \eqref{eq:eq:proof4}, and \eqref{eq:eq:proof5}, we obtain
\begin{align}\label{eq:eq:proof6}
	&\|x_{n+1}-x\|^2 + \|u_{n+1}-u\|^2 + \left(1-\frac{\tau}{2\beta}\right) \|v_{n+1}\|^2+2\tau_n\scal{Lv_{n+1}}{u_{n+1}-u} \nonumber\\
	& \leq  \|x_{n}-x\|^2 + \|u_{n}-u\|^2 - \|u_{n+1}-u_n\|^2 +2\tau_{n-1}\scal{Lv_n}{u_{n+1}-u}+\tau_n^2\|Dz_{n+1}-Dx_n\|^2.
\end{align}	
Moreover
\begin{align}\label{eq:eq:proof7}
	- \|u_{n+1}-u_n\|&^2 +2\tau_{n-1}\scal{Lv_n}{u_{n+1}-u} \nonumber\\
    &= 2\tau_{n-1}\scal{Lv_n}{u_{n}-u}  +\tau_{n-1}^2\|Lv_n\|^2 -\|u_{n}-u_{n+1}+\tau_{n-1}Lv_n\|^2\nonumber\\
    &\leq  2\tau_{n-1}\scal{Lv_n}{u_{n}-u}  +\tau^2\|L\|^2\|v_n\|^2 -\|u_{n}-u_{n+1}+\tau_{n-1}Lv_n\|^2.
\end{align}

Therefore, from \eqref{eq:eq:proof6}, \eqref{eq:eq:proof7},   \eqref{eq:desls1}, and \eqref{def:gamman}, we deduce that
\begin{align}	
\Gamma_{n+1}(x,u)
	&  \leq \|x_{n}-x\|^2 + \|u_{n}-u\|^2 +2\tau_{n-1}\scal{Lv_n}{u_{n}-u} +\tau^2\|L\|^2\|v_n\|^2 \nonumber\\
	&\hspace{4.5cm}-(1-\theta^2)\|u_n-u_{n+1}+\tau_{n-1}Lv_n\|^2\nonumber\\
    & = \Gamma_n(x,u)
    -\left(1-\frac{\tau}{2\beta}-\theta^2-\tau^2\|L\|^2\right)\|v_n\|^2 -(1-\theta^2)\|u_n-u_{n+1}+\tau_{n-1}Lv_n\|^2.\nonumber
\end{align}

\eqref{prop:des15} This follows directly from \eqref{eq:propdes1}, \eqref{eq:propdestau}, the fact that, for every $n \in \N$, $\tau_n \leq \tau$ and that
\begin{equation}\label{e:deschi}
    1-\frac{\tau}{2\beta}-\theta^2-\tau^2\|L\|^2 >0\quad  \Leftrightarrow\quad  \tau < \frac{4\beta(1-\theta^2)}{1+\sqrt{1+16\beta^2\|L\|^2(1-\theta^2)}}. 
\end{equation}
\eqref{prop:des2} Let $n\in\N$. It follows from \eqref{eq:CS}, \eqref{eq:YI}, and $\tau_{n-1}\le \tau$ that
\begin{equation*}
    2\tau_{n-1} \scal{Lv_n}{u_n-u} \geq -\left(1-\frac{\tau}{2\beta}-\theta^2\right)\|v_n\|^2-\frac{\tau^2\|L\|^2}{(1-\tau/(2\beta)-\theta^2)}\|u_{n}-u\|^2.
\end{equation*}
Hence, the result follows from \eqref{def:gamman} and \eqref{eq:propdestau}, noticing that,
\begin{align*}
	&\Gamma_{n}(x,u)\geq 
	\|x_{n}-x\|^2 + \left(1-	\frac{\tau^2\|L\|^2}{(1-\tau/(2\beta)-\theta^2)} \right)\|u_{n}-u\|^2. 
\end{align*}
\end{proof}

The next theorem guarantees weak convergence of Algorithm~\ref{algo:algo1} to a solution to Problem~\ref{pro:main}.
\begin{teo}\label{teo:main}
In the context of Problem~\ref{pro:main}, let $\theta\in\left]0,1\right[$ and
let $\tau_{-1}=\tau\in\left]0,\chi(\theta)\right[$, where $\chi(\theta)$ is
defined in \eqref{eq:propdestau}. Let
$(x_{-1},x_0,z_0,u_0)\in\H\times X\times\H\times\G$ and consider the sequence
$(x_n,z_n,v_n,u_n)_{n\in\N}$ generated by Algorithm~\ref{algo:algo1} with
$(\tau_n)_{n\in\N}$ defined by Line Search~\ref{line search1}. Suppose that,
for every $x\in X$, there exists $\rho>0$ such that
$D\big(\mathcal{B}(P_{\overline{\dom}A}x,\rho)\cap\dom A\big)$ is bounded.
Moreover, suppose that at least one of the following additional statements
holds:
\begin{enumerate}[label=(\alph*)]
\item\label{teo:hypa} $\displaystyle\liminf_{n\to\infty}\tau_n=\underline{\tau}>0$.
\item\label{teo:hypb} $D$ is locally uniformly continuous on $X$ relative to
$\dom A$, i.e., $Dx_k-Dy_k\to0$ for every sequence $(x_k)_{k\in\N}$ in $X$
and every sequence $(y_k)_{k\in\N}$ in $\dom A$ such that
$(x_k)_{k\in\N}$ converges weakly and $\|x_k-y_k\|\to0$.
\item\label{teo:hypc} $D$ is uniformly continuous on every weakly compact
subset of $\dom A\cup X$.
\end{enumerate}
Then $(x_n,u_n)_{n\in\N}$ converges weakly to a point in $\bm{Z}$.
\end{teo}
\begin{proof} For every $n\in\N$, set $\x_n=(x_n,z_n,v_n,u_n)$ and let $\tau_n$ be the
output of Line Search~\ref{line search1} at $\x_n$ with $\tilde\tau=\tau_{n-1}$.
Note that $x_0\in X$ and $x_{n+1}\in\ran P_X=X$ for every $n\in\N$.
Moreover, if $x_n\in\overline{\dom}A$, then assertion \ref{prop:LS1} of
Proposition~\ref{prop:LS} holds at $\x_n$. Otherwise, since $x_n\in X$, our
assumption yields assertion \ref{prop:LS2}. Therefore, by induction, Line
Search~\ref{line search1} terminates at every iteration, $(\tau_n)_{n\in\N}$
is a well-defined sequence in $\left]0,\tau\right]$, and \eqref{eq:algo1} and
\eqref{eq:defqppk} yield $u_{n+1}=\fu(\tau_n,\x_n)$ and
$z_{n+1}=\fz(\tau_n,\x_n)$, where $\fu$ and $\fz$ are computed with
$\tilde\tau=\tau_{n-1}$. In particular, \eqref{eq:line search1} with
$\alpha_k=\tau_n$ is \eqref{eq:desls1}. Moreover, in view of $\tau < \chi (\theta)$ and \eqref{e:deschi}, Proposition~\ref{prop:des}\eqref{prop:des1} and \cite[Lemma~3.1]{Comb2001} imply that, for every $(x,u)\in \bm{Z}$,
	$(\Gamma_n(x,u))_{n \in \N}$ is convergent,
		\begin{equation*}
			\sum_{n\in\N}\|u_{n+1}-u_n-\tau_{n-1}Lv_n\|^2<\pinf \quad\text{and}\quad \sum_{n\in\N}\|z_{n+1}-x_n\|^2=\sum_{n\in\N}\|v_{n+1}\|^2<\pinf.
		\end{equation*}In particular, from the continuity of $L$ we deduce
			\begin{equation}\label{eq:xuto0}
	v_{n+1}=	z_{n+1}-x_n \to 0 \quad \textnormal{ and } \quad  u_{n+1}-u_n  \to 0.
	\end{equation}
		Since $(\Gamma_{n}(x,u))_{n \in \N}$ is convergent, we conclude that $(x_n,u_n)_{n \in \N}$  is a bounded sequence in view of Proposition~\ref{prop:des}\eqref{prop:des2}. 
        Observe that, we deduce from \eqref{eq:algo1} that	
         \begin{equation}\label{eq:convtoZ}
			\begin{aligned}
				t_n:=	\frac{u_n-u_{n+1}}{\tau_n}-L\left(v_{n+1}-\frac{\tau_{n-1}}{\tau_n}v_n\right) &\in B^{-1}u_{n+1}-Lz_{n+1},\\
				s_n:=    \frac{x_n-z_{n+1}}{\tau_n}-(Cx_n-Cz_{n+1}+Dx_n-Dz_{n+1}) &\in (A +C+D)z_{n+1}+L^*u_{n+1}.
			\end{aligned}
		\end{equation} We now prove that, if \ref{teo:hypa}, \ref{teo:hypb}, or \ref{teo:hypc} holds, then every weak cluster point of $(x_n,u_n)_{n \in \N}$ is in $\bm{Z}$. 
         
\ref{teo:hypa} There exist $\varepsilon>0$ and $n_0\in\N$ such that
$\tau_n\ge\varepsilon$ for every $n\ge n_0$. In view of \eqref{eq:xuto0}, the Lipschitz continuity of $L$ and $C$, and \eqref{eq:desls1}, we conclude that $t_n \to 0$ and $s_n \to 0$. Moreover, since  $A+D$ is maximally monotone and $C$ is $\beta$-cocoercive, it follows from \cite[Corollary~25.5(i)]{bauschkebook2017} 
		that $A+C+D$ is maximally monotone. Furthermore, it follows  from 
		\cite[Proposition~2.7]{briceno2011SIAM} that the operator $\bm{A}\colon (y,v) \mapsto \big((A+C+D)y+L^*v\big)\times \big(B^{-1}v-Ly\big)  $ is maximally monotone. Therefore, the weak-strong closure of $\gra \bm{A}$ \cite[Proposition~20.38]{bauschkebook2017}, \eqref{eq:xuto0}, and \eqref{eq:convtoZ} imply that every weak cluster point of $(x_n,u_n)_{n \in \N}$ belongs to $\bm{Z}$.

\ref{teo:hypb} Let $(\overline{x},\overline{u})$ be a weak cluster point of
$(x_n,u_n)_{n\in\N}$, say $(x_{n_k},u_{n_k})\weak(\overline{x},\overline{u})$.
If $\liminf_{k\to\infty}\tau_{n_k}>0$, the argument of case \ref{teo:hypa}
applied to $(x_{n_k},u_{n_k})_{k\in\N}$ yields
$(\overline{x},\overline{u})\in\bm{Z}$. Otherwise, passing to a subsequence
if necessary, we assume that $\tau_{n_k}\to0$ and $\tau_{n_k}<\tau$ for every
$k\in\N$. Let $n\in\{n_k\}_{k\in\N}$, and define $\widetilde{\tau}_n=\tau_n/\gamma$ and
\begin{equation}\label{eq:defxtut}
\begin{aligned}
\widetilde{u}_{n+1}	&= J_{\widetilde{\tau}_n B^{-1}} \big(u_n+\widetilde{\tau}_n Lx_n+\tau_{n-1}Lv_n\big)\\
\widetilde{z}_{n+1}	&= J_{\widetilde{\tau}_n A} \big(x_n-\widetilde{\tau}_n (L^* \widetilde{u}_{n+1}+Cx_n+Dx_n)\big).\end{aligned}
\end{equation}
Since $\tau_n<\tau$, the step-size $\widetilde{\tau}_n$ was 
rejected by Line Search~\ref{line search1}, i.e.,
\begin{equation}\label{eq:lsfail} \widetilde{\tau}_n^2\|D\widetilde{z}_{n+1}-Dx_n\|^2> \theta^2 (\|\widetilde{z}_{n+1}-x_n\|^2 + \| u_n-\widetilde{u}_{n+1}+\tau_{n-1}Lv_n \|^2).
\end{equation} 
Since $ \tau_n < \widetilde{\tau}_n$,  it follows from Lemma~\ref{lemma}\eqref{lemmai} with $z=u_n+\tau_{n-1}Lv_n$ and $y=-Lx_n$, that
        		\begin{align*}
			\frac{\| u_n +\tau_{n-1}Lv_n  - \widetilde{u}_{n+1}\|}{\widetilde{\tau}_n} &\leq \frac{\|u_n  +\tau_{n-1}Lv_n - u_{n+1}\|}{\tau_n}\quad \\
            &\Leftrightarrow\quad  \gamma \| u_n  - \widetilde{u}_{n+1}+\tau_{n-1}Lv_n \| \leq \| u_n  - u_{n+1}+\tau_{n-1}Lv_n \|.
	\end{align*}	
        Therefore, \eqref{eq:xuto0} and the continuity of $L$ yield						\begin{equation}\label{eq:utto0}
			u_n - \widetilde{u}_{n+1} \to 0.
		\end{equation}		
        On the other hand, it follows from the nonexpansiveness of $J_{\widetilde{\tau}_n A}$ and \eqref{eq:defxtut} that
        \begin{align*}
        \| x_n - \widetilde{z}_{n+1}\|&\le \| x_n - J_{\widetilde{\tau}_n A }(x_n-\widetilde{\tau}_n(L^*u_{n+1}+Cx_n+Dx_n))\|\\
        &\hspace{4cm}+\|  J_{\widetilde{\tau}_nA }(x_n-\widetilde{\tau}_n(L^*u_{n+1}+Cx_n+Dx_n))-\widetilde{z}_{n+1}\|\\
        &\le \| x_n - J_{\widetilde{\tau}_n A }(x_n-\widetilde{\tau}_n(L^*u_{n+1}+Cx_n+Dx_n))\|+
        \widetilde{\tau}_n\|  L^*(u_{n+1}-\widetilde{u}_{n+1})\|.
        \end{align*}
        Hence, Lemma~\ref{lemma}\eqref{lemmai} with $z=x_n$ and $y=L^*u_{n+1}+Cx_n+Dx_n$ and $ \tau_n < \widetilde{\tau}_n$ imply
        \begin{align*}
			\frac{\|x_n - \widetilde{z}_{n+1}\|}{\widetilde{\tau}_n} &\leq \frac{\| x_n - z_{n+1}\|}{\tau_n}+\|  L^*(u_{n+1}-\widetilde{u}_{n+1})\|\\
             &\leq \frac{\| x_n - z_{n+1}\|}{\tau_n}+\|L\|(\|  u_{n+1}-u_{n}\|+\|  u_{n}-\widetilde{u}_{n+1}\|),
	\end{align*}    
    which yields
     $
			\gamma\| x_n - \widetilde{z}_{n+1}\|
             \leq \| x_n - z_{n+1}\|+\tau_{n}\|L\|(\|  u_{n+1}-u_{n}\|+\|  u_{n}-\widetilde{u}_{n+1}\|).
	$
    Therefore, it follows from \eqref{eq:xuto0} and \eqref{eq:utto0} that
    \begin{equation}\label{eq:xtto0}
			x_n - \widetilde{z}_{n+1} \to 0.
		\end{equation}In view of \eqref{eq:xtto0} and \eqref{eq:utto0}, we have
$(\widetilde{z}_{n_k+1},\widetilde{u}_{n_k+1})\weak(\overline{x},\overline{u})$.
Moreover, since $(x_{n_k})_{k\in\N}$ is a weakly convergent sequence in $X$,
$(\widetilde{z}_{n_k+1})_{k\in\N}$ is a sequence in $\dom A$, and
\eqref{eq:xtto0} holds, it follows from \ref{teo:hypb} that $D \widetilde{z}_{n_k+1} - D x_{n_k} \to 0$. Therefore, in view of \eqref{eq:lsfail}, we deduce 
			\begin{equation*}
			\frac{x_{n_k}-\widetilde{z}_{n_k+1}}{\widetilde{\tau}_{n_k}} \to 0 \quad \textnormal{ and } \quad \frac{u_{n_k}-\widetilde{u}_{n_k+1}+\tau_{n_k-1}Lv_{n_k}}{\widetilde{\tau}_{n_k}} \to 0.
		\end{equation*}	
		Furthermore, from \eqref{eq:defxtut}, for every $k \in \N$, we have, 
				\begin{equation*}
					\begin{aligned}
				\frac{u_{n_k}+\tau_{n_k-1}Lv_{n_k}-\widetilde{u}_{n_k+1}}{\widetilde{\tau}_{n_k}}+L(x_{n_k}-\widetilde{z}_{n_k+1}) &\in B^{-1}\widetilde{u}_{n_k+1}-L\widetilde{z}_{n_k+1}\\
				\frac{x_{n_k}-\widetilde{z}_{n_k+1}}{\widetilde{\tau}_{n_k}}-(Cx_{n_k}-C\widetilde{z}_{n_k+1}+Dx_{n_k}-D\widetilde{z}_{n_k+1}) &\in (A +C+D)\widetilde{z}_{n_k+1}+L^*\widetilde{u}_{n_k+1},
			\end{aligned}
		\end{equation*}
and, arguing as in case \ref{teo:hypa} with
$(\widetilde{z}_{n_k+1},\widetilde{u}_{n_k+1})\weak(\overline{x},\overline{u})$,
we conclude that $(\overline{x},\overline{u})\in\bm{Z}$.

\ref{teo:hypc} It is enough to show that \ref{teo:hypc} implies
\ref{teo:hypb}. Let $(x_k)_{k\in\N}$ be a sequence in $X$ and
$(y_k)_{k\in\N}$ be a sequence in $\dom A$ such that
$x_k\weak\overline{x}$ and $\|x_k-y_k\|\to0$. Since $X$ is closed and
convex, it is weakly closed and, hence, $\overline{x}\in X$. Therefore,
$\{\overline{x}\}\cup\bigcup_{k\in\N}\{x_k,y_k\}$ is a weakly compact
subset of $\dom A\cup X$, as the union of two sequences weakly converging
to $\overline{x}$ together with their limit. Thus, \ref{teo:hypc} yields
$Dx_k-Dy_k\to0$, and \ref{teo:hypb} holds.

        Finally, to conclude the proof using \cite[Lemma~2.46]{bauschkebook2017}, we will show that $(x_n,u_n)_{n \in \N}$ possesses at most one weak cluster point. Let $(x_{n_k},u_{n_k})_{k \in \N}$ and $(x_{n_m},u_{n_m})_{m \in \N}$ be two subsequences of $(x_n,u_n)_{n \in \N}$ such that 
        $(x_{n_k},u_{n_k})\weak 
(x,u) \in \H \times \G$ 
and $(x_{n_m},u_{n_m})\weak (y,v) \in \H \times \G$. In particular, since $(x,u) \in \bm{Z}$ and $(y,v) \in \bm{Z}$, the sequences $(\Gamma_n(x,u))_{n 
	\in 
	\N}$ and $(\Gamma_n(y,v))_{n 
	\in 
	\N}$ are convergent. Suppose that $\Gamma_n(x,u) \to \ell_1 \in \R$ and  $\Gamma_n(y,v) \to \ell_2 \in \R$. It follows from \eqref{def:gamman} that, for every $n \in \N$,
\begin{align}
	&\Gamma_n(x,u)-\Gamma_n(y,v) \nonumber\\
    =\,&
	\|x_{n}-x\|^2 + \|u_{n}-u\|^2        +2\tau_{n-1}\scal{Lv_n}{v-u}-\|x_{n}-y\|^2 -\|u_{n}-v\|^2 \nonumber\\
	=\,& \|x\|^2+
	2\scal{x_n}{y-x}-\|y\|^2+\|u\|^2
	+2\scal{u_n}{v-u}-\|v\|^2+2\tau_{n-1}\scal{Lv_n}{v-u}. \label{eq:ell1elle20}
\end{align}
In particular, taking the limit along $n_k\to { +}\infty$, since $v_n \to 0$, we deduce
\begin{align}\label{eq:ell1ell21}
	\ell_1-\ell_2&=\|x\|^2
    +2\scal{x}{y-x}-\|y\|^2+\|u\|^2+2\scal{u}{v-u}-\|v\|^2=-\|x-y\|^2-\|u-v\|^2.
\end{align}
Analogously, taking the limit along $n_m\to{ +}\infty$ in \eqref{eq:ell1elle20}, we obtain
\begin{equation}\label{eq:ell1ell22}
	\ell_1-\ell_2=\|x-y\|^2+\|u-v\|^2.
\end{equation}
Combining \eqref{eq:ell1ell21} and \eqref{eq:ell1ell22}, we have that $\|x-y\|^2+\|u-v\|^2 = 0$, which implies $x=y$ and $u=v$. This completes the proof. 
 \end{proof}
\begin{rem}\label{rem:convergence} 
\begin{enumerate}
\item
The first assumption in Theorem~\ref{teo:main} holds automatically at every
$x\in X\cap\overline{\dom}A$, by the continuity of $D$ on $\dom A\cup X$.
In particular, it holds if $X\subset\overline{\dom}A$, and also under any of
the conditions in Remark~\ref{rem:LS}.

\item If $D$ is Lipschitz continuous on bounded subsets of $\dom A\cup X$,
then the first assumption holds by Remark~\ref{rem:LS} and \ref{teo:hypa} holds. Indeed, since $(x_n,u_n,v_n)_{n\in\N}$
is bounded, the nonexpansiveness of the resolvents implies that the iterates
and all the trial points of Line Search~\ref{line search1} lie in a bounded
subset of $\dom A\cup X$, on which $D$ is $\ell$-Lipschitz for some $\ell>0$.
Hence every $\alpha_k\le\theta/\ell$ is accepted and
$\tau_n\ge\min\{\tau,\gamma\theta/\ell\}$. In finite dimensions, it suffices that
every point of $\overline{\dom A\cup X}$ has a neighborhood on which
$D|_{\dom A\cup X}$ is Lipschitz continuous.

\item Assumption \ref{teo:hypc} holds if $D$ is uniformly continuous on
every weakly compact subset of $\overline{\conv}(\dom A\cup X)$, as assumed
in \cite[Theorem~3.4]{BricenoRoldan4op}. In finite dimensions,
\ref{teo:hypc} follows from the continuity of $D$ on $\dom A\cup X$.
\item When $\G=\{0\}$, Algorithm~\ref{algo:algo1} with
Line Search~\ref{line search1} reduces to the method in
\cite[Theorem~2.3]{BricenoDavis2018} with line search, and
Theorem~\ref{teo:main} recovers its convergence under \ref{teo:hypa},
without assuming $X\subset\dom A$ or the maximal monotonicity of $D$.
Under \ref{teo:hypc}, uniform continuity is required on weakly compact
subsets of $\dom A\cup X$ instead of $X$, which is consistent with the proof
of \cite[Theorem~2.3]{BricenoDavis2018}, where $D$ is also evaluated at
points of $\dom A$.
\item If, in addition, $C=0$, we recover the method in
\cite{Tseng2000SIAM} with the step-size rule \cite[(2.4)]{Tseng2000SIAM}.
In this case, Theorem~\ref{teo:main} recovers
\cite[Theorem~3.4(a)--(b)]{Tseng2000SIAM}, since \ref{teo:hypa} and
\ref{teo:hypb} are the conditions assumed there, and it does not need the
local boundedness of $x\mapsto\min_{w\in(A+D)x}\|w\|$, the inclusion
$X\subset\dom A$, or the maximal monotonicity of $D$.

\item\label{rem:convergencesepsize} Algorithm~\ref{algo:algo1} employs the same step-size for the primal and dual
variables. Primal-dual algorithms usually allow for independent primal and dual step-sizes \cite{ChambollePock2011,Condat13,Roldan20254op,Vu13}. In our setting, we can rescale the dual update to induce a fixed ratio between the primal and dual step-sizes. More precisely, let $\omega\in\left]0,+\infty\right[$
and define
\begin{equation*}
\widetilde{L}=\sqrt{\omega} L
\quad\text{and}\quad
\widetilde{B}\colon\G\to 2^{\G}\colon y\mapsto \frac{1}{\sqrt{\omega}}B\left(\frac{y}{\sqrt{\omega}}\right).
\end{equation*}
Note that $\widetilde{L}$ is linear and bounded with $\|\widetilde{L}\|=\sqrt{\omega}\|L\|$,
$\widetilde{B}$ is maximally monotone, and
$\widetilde{B}^{-1}\colon u \mapsto\sqrt{\omega}B^{-1}\left(\sqrt{\omega}u \right)$.
Hence, for every $(x,u)\in\H\times\G$,
\begin{equation*}
\begin{cases}
0\in Ax+Cx+Dx+L^*u\\
0\in B^{-1}u-Lx
\end{cases}
\quad\Leftrightarrow\quad
\begin{cases}
0\in Ax+Cx+Dx+\widetilde{L}^*(u/\sqrt{\omega})\\
0\in \widetilde{B}^{-1}(u/\sqrt{\omega})-\widetilde{L}x,
\end{cases}
\end{equation*}
that is, $(x,u)\in \bm{Z}$ if and only if $(x,u/\sqrt{\omega})$ solves Problem~\ref{pro:main}
with $L$ and $B$ replaced by $\widetilde{L}$ and $\widetilde{B}$.
Moreover, by \cite[Corollary~23.26]{bauschkebook2017}, we have
$J_{\widetilde{B}^{-1}}=\frac{1}{\sqrt{\omega}}J_{\omega B^{-1}}\circ(\sqrt{\omega}\id )$.
Therefore, the sequence $(v_{n-1},\tilde u_n,z_n,x_n)_{n \in \N}$ generated by Algorithm~\ref{algo:algo1} applied to the
rescaled problem can be written equivalently, by setting $u_n=\sqrt{\omega}\widetilde{u}_n$ and $\sigma_n=\omega \tau_n$, as
\begin{equation}\label{e:alg_gamma}
(\forall n\in\mathbb{N})\quad
\begin{cases}
v_n=z_n-x_{n-1}\\
u_{n+1}=J_{\sigma_nB^{-1}}\big(u_n+\sigma_nLx_n+\sigma_{n-1}Lv_n\big)\\
z_{n+1}=J_{\tau_nA}\big(x_n-\tau_n(L^*u_{n+1}+Cx_n+Dx_n)\big)\\
x_{n+1}=P_X\big(z_{n+1}-\tau_n(Dz_{n+1}-Dx_n)\big).
\end{cases}
\end{equation}
In this scenario, $\tau_n$ is determined by Line Search~\ref{line search1}, where
$q_k=J_{{\omega}\alpha_kB^{-1}}\big(u_n+{\omega}(\alpha_kLx_n+\tau_{n-1}Lv_n)\big)$ and  \eqref{eq:line search1} reduces to
\begin{equation}\label{eq:linsearchrem}
\alpha_k^2\|Dp_k-Dx_n\|^2\leq\theta^2\Big(\|p_k-x_n\|^2
+{\frac{1}{\omega}}\big\|u_n-q_k+{\omega}\tau_{n-1}Lv_n\big\|^2\Big).
\end{equation}
Since $A$, $C$, $D$, and $X$ are not modified, the assumptions of
Theorem~\ref{teo:main} hold for the rescaled problem, and we conclude that
$(x_n,u_n/\sqrt{\omega})_{n\in\mathbb{N}}$ converges weakly to a solution of the rescaled
problem, i.e., $(x_n,u_n)_{n\in\mathbb{N}}$ converges weakly to a point in $\bm{Z}$.  In the case where $X=\H$ and $D$ is $\zeta$-Lipschitz continuous, we have that, for every $n \in \N$, \eqref{eq:linsearchrem} holds for $\alpha_k = \tau$ and $\theta = \zeta \tau$. In this case, the recurrence in \eqref{e:alg_gamma} coincides with the method proposed in \cite[Theorem~4.2]{Roldan20254op}  
and \eqref{eq:propdestau}  coincides with the condition guaranteeing the convergence of this method, in which the parameter $\varepsilon$ of \cite[Theorem~4.2]{Roldan20254op} is taken as $\tau/(2\beta)$. As a consequence, we recover the convergence of the Condat--V\~u \cite{Condat13,Vu13}, Chambolle--Pock \cite{ChambollePock2011}, FBHF, and forward-backward-forward methods \cite{Tseng2000SIAM}. Moreover, since $\|\widetilde L\|=\sqrt{\omega}\|L\|$, condition \eqref{eq:propdestau} for the rescaled problem reads
$\tau\in\left]0,\chi_\omega(\theta)\right[$, where
$$\chi_\omega(\theta)=\frac{4\beta(1-\theta^2)}{1+\sqrt{1+16\beta^2\omega\|L\|^2(1-\theta^2)}}.$$
\item Note that every step-size accepted by the purely primal test
$\alpha_k\|Dp_k-Dx_n\|\le\theta\|p_k-x_n\|$ also satisfies \eqref{eq:line search1}. Hence, if the step-sizes are generated by this test, Proposition~\ref{prop:des} and case~\ref{teo:hypa} of Theorem~\ref{teo:main} remain valid. However, the dual term in \eqref{eq:line search1} is crucial in case~\ref{teo:hypb} of  Theorem~\ref{teo:main}, where $\tau_n\to0$ and the proof uses that the rejected step-size $\tau_n/\gamma$
violates \eqref{eq:line search1} (see \eqref{eq:lsfail}). Indeed, as we will see in Example~\ref{ex:primalLSconvL}, it may happen that the primal line search terminates at every iteration, but $(x_n,u_n)_{n\in\N}$ converges to a point which is not in $\bm{Z}$.
\end{enumerate}
\end{rem}

\begin{exam}\label{ex:primalLSconvL} In the context of Problem~\ref{pro:main}, set $\H=\R^2$, $\G=\R$, $X=\R^2$,
$C=0$, $B=\id$, $
	A=N_{\R\times\RP}$, $
	L\colon(t,r)\mapsto r$, and $
	D\colon(t,r)\mapsto\big(\operatorname{sgn}(t)\sqrt{|t|},0\big)$. The assumptions of Theorem~\ref{teo:main} hold and $\bm{Z}=\{((0,0),0)\}$ . Let $\theta=3/4$, $\tau=\tau_{-1}=15/64$, $\beta=1$, and $\gamma = 7/8$. For this parameters
\eqref{eq:propdestau} holds and Theorem~\ref{teo:main}\ref{teo:hypc}
guarantees that the sequence generated by Algorithm~\ref{algo:algo1} with
Line Search~\ref{line search1} converges to $((0,0),0)$. However, if we replace \eqref{eq:line search1} by a purely primal test, and we consider $x_0=(1/16,0)$, $u_0=1$, and $z_0=x_{-1}$, it can be proved that $x_n \to (0,0)$ and $u_n \to u^* >0$ (see the Appendix). Hence, $(x_n,u_n)\to((0,0),u^*)\notin\bm{Z}$. Moreover, in terms of the proof of
Theorem~\ref{teo:main}, with $\widetilde\tau_n=\tau_n/\gamma$, we have
\begin{equation*}
	\frac{u_n-\widetilde u_{n+1}+\tau_{n-1}Lv_n}{\widetilde\tau_n}
	\to u^*,
	\qquad\text{and}\qquad
	\frac{x_n-\widetilde z_{n+1}}{\widetilde\tau_n}\to0.
\end{equation*}
That is, without the dual term, the failure of the test for
$\widetilde\tau_n$ only controls the primal residual, and the summability of
the step-sizes, dictated by the singularity of $D$ at the primal solution,
freezes the dual sequence.
\end{exam}

\section{Application to convex optimization problems}\label{sec:aplication}
In this section we apply our proposed algorithm for solving 
convex problems involving nonlinear composite terms. These optimization problems encompass applications in data science \cite{CombettesPesquet2021strategies}, machine learning \cite{Nocedal2018}, image processing
\cite{Berteroetal2010,BotHendrich2014TV,Briceno2011ImRe,BurgerSawatzkySteidl2016,chambolle2016AN,TeuberSteidlChan2013}, among others.

Let $\mathcal K$ be a real Hilbert space, and let $K\subset\mathcal K$ be
a nonempty closed convex cone. A function $\phi\colon\H\to\mathcal K$ is \emph{$K$-convex} if
\[
(\forall (x,y)\in\H^2)(\forall\lambda\in[0,1])\quad
\lambda \phi(x)+(1-\lambda)\phi(y)-\phi\big(\lambda x+(1-\lambda)y\big)\in K.
\]  Furthermore, $\psi\colon\mathcal K\to\RX$ is \emph{$K$-nondecreasing} if, for every $(x,y)\in\mathcal K^2$, $y-x\in K$ implies $\psi(x)\le\psi(y)$.

Consider the following convex optimization problem involving nonlinear compositions.
\begin{pro}\label{pro:opti_nonlinear}
	Let $\H$, $\G$, and $\mathcal{K}$ be real Hilbert spaces, and let $K \subset \mathcal{K}$ be a closed convex cone. Let $f \in \Gamma_0(\H)$, $g \in \Gamma_0(\G)$, and $\psi \in \Gamma_0(\mathcal{K})$ be such that $\psi$ is $K$-nondecreasing. Let $h\colon \H \to \R$ be a differentiable convex function with a $\beta^{-1}$-Lipschitz gradient for some $\beta\in \RPP$. Let $L \colon \H \to \G$ be a bounded linear operator, and let $d\colon \H \to \mathcal{K}$ be a continuously differentiable $K$-convex function. The problem is 
\begin{equation} \label{eq:prob_nonlinear}
    \min_{x \in \H} f(x) + g(Lx) + h(x) + \psi(d(x)),
\end{equation}
	under the assumption that its solution set is nonempty.
\end{pro}
\begin{rem} \label{r:K} Some comments regarding the assumptions in Problem~\ref{pro:opti_nonlinear} are in order.
\begin{enumerate}
    \item Since $\psi \in \Gamma_0(\mathcal{K})$ is $K$-nondecreasing and $d$ is $K$-convex, $\psi \circ d$ is convex \cite[Theorem 2.1.3(vi)]{zalinescu2002convex}, and  problem \eqref{eq:prob_nonlinear} is convex.
    \item \label{r:K-2} Let $K^*:= \{s \in \mathcal{K}^*: \scal{x}{s}\geq 0 \quad\forall x \in K\}$ denote the dual cone of $K$. It holds that $\psi \in \Gamma_0(\mathcal{K})$ is $K$-nondecreasing if and only if $\dom \psi^* \subseteq K^*$ (see \cite[Lemma 2.8]{MR4888251}). Indeed, suppose $\psi$ is $K$-nondecreasing and let $z \notin K^*$. There exists $\hat y\in K$ with $\scal{z}{\hat y}<0$. Fix $y_0\in\dom\psi$,
which exists because $\psi$ is proper. For every $t>0$, $y_0-(y_0-t\hat y)=t\hat y\in K$, so
$\psi(y_0-t\hat y)\le\psi(y_0)<\pinf$. Therefore
\begin{equation*}
\psi^*(z)\;\ge\;\scal{z}{y_0-t\hat y}-\psi(y_0-t\hat y)\;\ge\;\scal{z}{y_0}-t\,\scal{z}{\hat y}-\psi(y_0)
\;\to\;\pinf \text{ as } t\to\pinf,
\end{equation*}
then $z\notin\dom\psi^*$.   Conversely, suppose $\dom \psi^* \subseteq K^*$. Let $x_2 - x_1 \in K$ and $z \in \dom \psi^*$. Then, in view of the definition of the dual cone, \begin{equation*}
    \psi^{**}(x_2) \geq \scal{x_2}{z} - \psi^*(z)  = \scal{x_2-x_1}{z} + \scal{x_1}{z} - \psi^*(z) \geq \scal{x_1}{z} - \psi^*(z).
\end{equation*} Taking the supremum over $z \in \dom \psi^*$ and using that $\psi = \psi^{**}$, it follows that $\psi(x_2) \geq \psi(x_1)$, which proves the claim.
\end{enumerate}
\end{rem}
This problem can be solved by the methods in \cite{BricenoDavis2018,BricenoRoldan4op,Tseng2000SIAM} using an equivalent primal--dual formulation and treating the linear term as a (Lipschitz) continuous operator. Instead, based on Algorithm~\ref{algo:algo1}, we propose the following algorithm to solve it. This algorithm takes advantage of the linearity, producing numerical advantages as we will see in Section~\ref{sec:num}.
\begin{algo}\label{algo:algononlinear}
	In the context of Problem~\ref{pro:opti_nonlinear}, let $(x_{-1},x_0,z_0) \in \H^3$, $w_0\in \mathcal{K}$, and $u_0 \in \G$. Let $\tau\in\RPP$, let $\tau_{-1}\in\left]0,\tau\right]$, and let $(\theta,\gamma) \in ]0,1[^2$. Set $v_0 = z_0-x_{-1}$ and consider the following recurrence:
\begin{equation*}
    (\forall n\in\N)\quad 
    \left\lfloor 
    \begin{array}{l} 
        \left\lfloor  
\begin{array}{l} 
    \textbf{For } k = 0, 1, \dots \textbf{ define } \\[1ex]
    \alpha_k = \tau \gamma^k \\[1ex]
    q_k = \prox_{\alpha_k g^{*}} \big(u_n + \alpha_k Lx_n + \tau_{n-1}Lv_n\big)\\[1ex]
    p_k = \prox_{\alpha_k f} \big(x_n - \alpha_k (L^* q_k + \nabla h(x_n) + [\nabla d(x_n)]^* w_n)\big)\\[1ex]
    r_k = \prox_{\alpha_k \psi^*} \big(w_n + \alpha_k d(x_n)\big)\\[1ex]
    \textbf{if } \alpha_k^2\big(\|[\nabla d(p_k)]^* r_k - [\nabla d(x_n)]^* w_n\|^2 + \|d(p_k) - d(x_n)\|^2\big) \\
    \hspace{1cm} \leq \theta^2 \big(\|p_k - x_n\|^2 + \|r_k - w_n\|^2 + \| u_n - q_k + \tau_{n-1}Lv_n \|^2\big) \\[1ex]
    \quad \textbf{break} \textnormal{ and return } \alpha_k,\ q_k,\ p_k,\ r_k
\end{array}
        \right.\\[4ex]
        \tau_n = \alpha_k \\[1ex]
        (u_{n+1},z_{n+1},y_{n+1}) = (q_k,p_k,r_k)\\[1ex]
        x_{n+1} = P_{\overline{\dom} f}\big(z_{n+1} - \tau_n\big([\nabla d(z_{n+1})]^* y_{n+1} - [\nabla d(x_n)]^* w_n\big)\big)\\[1ex]	
        w_{n+1} = P_{\overline{\dom} \psi^*}\big(y_{n+1} + \tau_n\big(d(z_{n+1}) - d(x_n)\big)\big)\\[1ex]
        v_{n+1} = z_{n+1} - x_{n}.
    \end{array}
    \right.
\end{equation*}
\end{algo}
The following theorem establishes the convergence of Algorithm~\ref{algo:algononlinear}.
\begin{teo}\label{teo:conv_nonlinear}
	In the context of Problem~\ref{pro:opti_nonlinear}, let $\theta\in~]0,1[$ and let $\tau_{-1} = \tau \in ]0,\chi(\theta)[$ where $\chi(\theta)$ is defined in \eqref{eq:propdestau}.
    Let $(x_{-1},x_0,z_0) \in \H^3$, $w_0 \in \mathcal{K}$, and $u_0 \in \G$. Consider the sequence $(v_n,z_n,x_n,u_n, y_n, w_n)_{n \in \N}$ generated by Algorithm~\ref{algo:algononlinear}. Suppose that the following qualification conditions hold:
\begin{equation} \label{eq:CQ_nonlinear_1}
    0 \in \sri(L(\dom f) - \dom g),
\end{equation}
and
\begin{equation} \label{eq:CQ_nonlinear_2}
    0 \in \sri\big(d(\dom f \cap L^{-1}(\dom g)) - \dom \psi + K\big).
\end{equation}
Moreover, assume that at least one of the following statements is satisfied:
		\begin{enumerate}[label=(\alph*)]
			\item 
            $\displaystyle \liminf_{n \to \infty} \tau_n = \underline{\tau} >0$.
			\item \label{teo:hypb2} The operator $\nabla d$ is uniformly continuous on every
weakly compact subset of $\overline{\dom}f$. \end{enumerate}
		Then, $(x_n)_{n \in \N}$ converges weakly to a solution to Problem~\ref{pro:opti_nonlinear}, and $(w_n)_{n \in \N}$ converges weakly to a dual solution of the nonlinear composition.
\end{teo}

\begin{proof} 
Since $\dom\psi^*\subseteq K^*$ in view of Remark~\ref{r:K}\eqref{r:K-2} and $d$ is $K$-convex, the
function $x\mapsto\scal{d(x)}{w}$ is convex for every $w\in\dom\psi^*$ \cite[Proposition~3.7]{Gissler2023}. Moreover, since $d$ is
continuously differentiable, we have  $  \nabla\big(\scal{d(\cdot)}{w}\big)(x)=\nabla d(x)^*w$. 

Under
\eqref{eq:CQ_nonlinear_1},  \cite[Theorem~16.47(i)]{bauschkebook2017} yields, for every $w\in\dom\psi^*$,
\begin{align*}
&\partial\big(
f+g\circ L+h+\scal{d(\cdot)}{w}
\big)(x) =
\partial f(x)+L^*\partial g(Lx)+\nabla h(x)
+[\nabla d(x)]^*w.
\end{align*}
On the other hand, \eqref{eq:CQ_nonlinear_2} and
\cite[Proposition~4.11]{Combari1994} imply that
\begin{align*}
&\partial\big(f+g\circ L+h+\psi\circ d\big)(x) =
\bigcup_{w\in\partial\psi(d(x))}
\partial\big(
f+g\circ L+h+\scal{d(\cdot)}{w}
\big)(x).
\end{align*}
Therefore, $\bar x$ solves
Problem~\ref{pro:opti_nonlinear} if and only if there exists
$\bar w\in\mathcal K$ such that
\begin{equation}\label{eq:KT_nonlinear_full}
\begin{cases}
0\in
\partial f(\bar x)+L^* \partial g(L\bar x)+\nabla h(\bar x)
+[\nabla d(\bar x)]^*\bar w,\\
\bar w\in\partial\psi(d(\bar x)).
\end{cases}
\end{equation}
Since $
\bar w\in\partial\psi(d(\bar x))$ is equivalent to $0\in\partial\psi^*(\bar w)-d(\bar x)$,
\eqref{eq:KT_nonlinear_full} can be equivalently written as
\begin{equation*}
\begin{cases}
0\in
\partial f(\bar x)+L^* \partial g(L\bar x)+\nabla h(\bar x)
+[\nabla d(\bar x)]^*\bar w,\\
0\in\partial\psi^*(\bar w)-d(\bar x).
\end{cases}
\end{equation*}

Now, set $\bm{\H} = \H \times \mathcal{K}$, $\bm X
=
\overline{\dom f}\times\overline{\dom\psi^*}$, and the operators $\bm{A},\bm{B}, \bm{L}, \bm{C}$, and $\bm{D}$ as follows:
\begin{equation}\label{eq:non-lin-operators}
    \begin{aligned}
        &\bm{A} \colon \HH \to 2^{\HH} \colon (x,w) \mapsto \begin{pmatrix} \partial f(x) \\ \partial \psi^*(w) \end{pmatrix},\\
        &\bm{B} \colon \G \to 2^{\G} \colon u \mapsto \partial g (u),\\
        &\bm{L} \colon \HH \to \G \colon (x,w) \mapsto Lx,\\
        &\bm{C} \colon \HH \to \HH \colon (x,w) \mapsto \begin{pmatrix} \nabla h(x) \\ 0 \end{pmatrix},\\
        &\bm{D} \colon \HH \to \HH \colon (x,w) \mapsto \begin{pmatrix} [\nabla d(x)]^*w \\ -d(x) \end{pmatrix}.
    \end{aligned}
\end{equation}
Therefore, Problem~\ref{pro:opti_nonlinear} is equivalent to find $\x \in \HH$ such that $0 \in (\bm{A} + \bm{C} + \bm{D} + \bm{L}^* \bm{B} \bm{L})\bm{x}$. Since $f \in \Gamma_0(\H)$ and $\psi^* \in \Gamma_0(\mathcal{K})$, the operator $\bm{A}$ is maximally monotone. Similarly, $\bm{B}$ is maximally monotone. We have that $\nabla h$ is $\beta^{-1}$-Lipschitz, thus, $\bm{C}$ is $\beta$-cocoercive on $\bm{\H}$. In addition, $\bm A+\bm D$ is the saddle subdifferential operator
associated with the closed proper saddle function \cite[Corollary~1]{Rockafellar70} \begin{equation*}
\mathcal L(x,w):=
\begin{cases}
f(x)+\scal{d(x)}{w}-\psi^*(w), & x\in\dom f,\ w\in\dom\psi^*,\\
\pinf, & x\notin\dom f,\\
\minf, & x\in\dom f,\ w\notin\dom\psi^*.
\end{cases}
\end{equation*}
Hence, $\bm A+\bm D$ is maximally monotone by Lemma~\ref{l:A+D}\ref{it:saddle}.  Note that $\bm{D}$ is continuous and, if \ref{teo:hypb2} holds,  $\bm{D}$  is uniformly continuous on every
weakly compact subset of $\overline{\dom} f \times \overline{\dom} \psi^*$. Moreover, since $\dom\partial f$ and $\dom\partial\psi^*$ are dense in $\dom f$ and $\dom\psi^*$, respectively, we have $\bm X\subset\overline{\dom}\bm A$ and, hence, the first assumption in Theorem~\ref{teo:main} holds (see Remark~\ref{rem:convergence}).
Finally, note that $J_{\tau \bm{A}}(x,w) = (\prox_{\tau f}x,\prox_{\tau \psi^*}w)$, $J_{\tau \bm B^{-1}} u = \prox_{\tau g^*} u$, and $P_{\bm X} (x,w)= (P_{\overline{\dom  }f }x,P_{\overline{\dom  }\psi^*}w)$. Hence, Algorithm~\ref{algo:algononlinear} corresponds to Algorithm~\ref{algo:algo1} in this setting and the result follows  from Theorem~\ref{teo:main}.
 \end{proof}
\begin{rem}
    The operator $\bm D$ in \eqref{eq:non-lin-operators} is monotone on $\H\times K^*$, because for $(w,w')\in (K^*)^2$, \begin{align*}
        & \scal{\bm D(x,w)-\bm D(x',w')}{(x-x',w-w')} \\
        = \, & \scal{w}{d(x')-d(x)-\nabla d(x)(x'-x)}  +\scal{w'}{d(x)-d(x')-\nabla d(x')(x-x')}\ge0
    \end{align*}
by $K$-convexity. 
However, the operator $\bm D$ may not be monotone over the whole space. Indeed, take $\H=\K=\R$, $K=\left[0,+\infty\right[$, and $d(x)=x^2$ for all $x\in\R$. Then,  $\bm D(x,w)=(2xw,-x^2)$ satisfies $\scal{\bm D(1,-1)-\bm D(0,-1)}{(1,0)}=-2<0$. Nevertheless, as shown above, the operator $\bm A+\bm D$ is still maximally monotone.
\end{rem}

\section{Numerical Experiments in Computed Tomography with Poisson Noise}\label{sec:num}

 In computed tomography (CT), the goal is to obtain a 2D or 3D representation of a target sample from a set of measured attenuation images, or \emph{projections}, under different views of such sample. From the acquired projections, a discrete representation of the sample is reconstructed usually by means of closed-form or iterative numerical strategies \cite{KakSlaney1988}. In parallel-beam CT, the X-rays generated at an external source travel in \emph{parallel} with a detector perpendicular to the rays trajectory with finite $m$ detecting spots, either arranged in a 2D plane for 3D imaging or in a 1D linear detector for the 2D setting. We will perform our numerical experiments related to the latter, with a circular trajectory of the source and detector relative to the sample within $n_\theta$ uniformly-spaced steps. We will also consider \emph{transmission} CT, where we are interested in recovering the \emph{attenuation coefficient} of the sample since the detectors measure the transmitted beam after its interaction with the sample.

Let $x\in\RP^N$ be the image to be reconstructed (with $N=n^2$ pixels) of attenuation coefficients of the sample, and let $b\in\mathbb{N}^M$ $(M = mn_\theta)$ be the measured data, considered here as the number of photons reaching the detector after sample interaction. Let $W \colon \RP^N \to \RP^M$ be the discrete linear Radon transform mapping $x$ to the line integrals of $x$ along each ray direction, and let
$I_0\in\RPP$ be the expected number of incident photons per detecting pixel. According to the Beer--Lambert law \cite{KakSlaney1988}, the expected number of detected photons decays as the mapping
\begin{equation}\label{eq:beerlambert}
	\mu(x) = I_0\,e^{-Wx}. 
\end{equation}

Since the measurements are photon counts, the data $b
$ are modeled as realizations of independent Poisson random variables per pixel with means $\mu(\overline{x})$, where $\overline{x}$ is the true image. Therefore, we are interested in the statistical inverse problem
\begin{equation*}
\mathbf{b} \sim \text{Poisson}\left(\mu(\overline{x})\right), 
\end{equation*}
where $b$ is a realization of the random vector $\mathbf{b}$.

In the Bayesian inversion framework, the negative log-likelihood of the data fidelity term is, up to an additive constant, the generalized Kullback--Leibler divergence
\begin{equation}\label{eq:KLct}
	\mathrm{KL}(b\,|\,\mu(x))=\sum_{i=1}^M\Big(b_i\ln\frac{b_i}{\mu_i(x)}-b_i+\mu_i(x)\Big)
	=\sum_{i=1}^M\Big(I_0\,e^{-(Wx)_i}+b_i(Wx)_i\Big)+\eta,
\end{equation}
where $\eta \in \R$ and with the convention $0\ln 0=0$. In our setting, we consider an underdetermined inverse problem with $M<N$, then the minimization of \eqref{eq:KLct} for $x$ over $\RP^N$ may overfit the noise \cite{BerteroBoccacci1998,Berteroetal2009}. Therefore, a regularization term is also considered and a usual choice is to  minimize $\lambda\,\mathrm{TV}(x)+\mathrm{KL}(b\,|\,\mu(x))$ over $\RP^N$, where $\mathrm{TV}$ is the isotropic total variation \cite{RudinOsherFatemi1992} and $\lambda\in\RPP$ is a regularization parameter, whose choice is a delicate issue. In the Bayesian setting, the TV term corresponds to assuming a Laplace distribution prior on the discrete gradient of the reconstructed image. 

For Poisson data, \cite{Berteroetal2010} shows that, for large counts, the expected value of
$\mathrm{KL}(b\,|\,\mu(\overline{x}))$ is approximately $M/2$, and proposes a discrepancy principle
for selecting $\lambda$, namely, $\lambda$ is chosen such that the corresponding solution $x_\lambda$
satisfies $\mathrm{KL}(b\,|\,\mu(x_\lambda))=\eta$, with $\eta=M/2$. In
\cite{TeuberSteidlChan2013}, the authors consider the equivalent constrained formulation, in which a
seminorm (e.g., $\mathrm{TV}$) is minimized subject to an I-divergence (i.e., Kullback--Leibler)
constraint, and the regularization parameter is recovered as the inverse of the Lagrange
multiplier associated to the constraint. Following this approach, we consider the problem
\begin{equation}\label{eq:ct_problem}
	\min_{\substack{x\in\RP^N\\ \mathrm{KL}(b\,|\,\mu(x))\le\eta}}\ \mathrm{TV}(x),
	\qquad\text{where}\quad \mathrm{TV}(x)=\|\nabla x\|_{1,2}=\sum_{j=1}^N|(\nabla x)_j|_2,
\end{equation}
$\nabla\colon\R^N\to\R^{N\times 2}$ is the discrete gradient defined by forward differences with
Neumann boundary conditions, and $|\cdot|_2$ is the Euclidean norm in $\R^2$. In contrast with the
models considered in \cite{Berteroetal2010,TeuberSteidlChan2013}, in which the expected counts depend
linearly on the image, in transmission CT the map $x\mapsto\mu(x)$ in \eqref{eq:beerlambert} is
nonlinear. Hence, the constraint in \eqref{eq:ct_problem} is a nonlinear composite term whose
projection is not available, and the methods in \cite{TeuberSteidlChan2013} cannot be applied
directly.

The optimization problem in \eqref{eq:ct_problem} is a
particular instance of Problem~\ref{pro:opti_nonlinear} for $f=\iota_{\RP^N}$, $g=\|\cdot\|_{1,2}$,
$L=\nabla$, $h=0$, $d\colon x\mapsto\mathrm{KL}(b\,|\,\mu(x))-\eta$,
$\psi=\iota_{\left]-\infty,0\right]}$, and $K=\RP$. In this setting, $\psi^*=\iota_{\RP}$, $d$ is
convex (and hence $K$-convex), since it is the composition of the convex function
$p\mapsto\sum_{i=1}^M(I_0e^{-p_i}+b_ip_i)$ with $W$, and
\begin{equation*}
	\nabla d(x)=W^\top\big(b-I_0\,e^{-Wx}\big),
\end{equation*}
where the exponential is taken componentwise. Note that each evaluation of $d$ and $\nabla d$
at a new point requires one forward projection ($Wx$) and one back projection (the product with
$W^\top$), which are the most expensive operations in this problem. Since $h=0$, its gradient is
$\beta^{-1}$-Lipschitz continuous for every $\beta\in\RPP$, and we set $\beta=+\infty$ (i.e.,
$1/\beta=0$) in all the step-size conditions. The operator $\nabla d$ is not Lipschitz continuous on
$\R^N$, but it is Lipschitz continuous on $\RP^N$ with constant $I_0\|W\|^2$, since $W\ge0$.
However, this constant is very conservative and, more importantly, the operator
$(x,w)\mapsto(w\nabla d(x),-d(x))$, which appears in the primal-dual formulation of
\eqref{eq:ct_problem}, is neither the gradient of a function nor globally Lipschitz continuous on $\RP^N\times\RP$, because of the product
$w\nabla d(x)$ and since a bound on the multiplier $w$ is not known a priori. Hence, methods
involving line searches are appropriate for this problem. Moreover, since $\dom g=\R^{N\times2}$,
\eqref{eq:CQ_nonlinear_1} holds and \eqref{eq:CQ_nonlinear_2} reduces to the Slater condition,
i.e., the existence of $x\in\RP^N$ such that $\mathrm{KL}(b\,|\,\mu(x))<\eta$, which holds when the
data can be fitted below the noise level. In addition, since $\nabla d$ is Lipschitz continuous on
$\dom f=\RP^N$, it is uniformly continuous on $\dom f$. Hence, Theorem~\ref{teo:conv_nonlinear}
guarantees the convergence of Algorithm~\ref{algo:algononlinear} in this case.

In order to balance the scales of the primal and dual variables, we multiply the objective function
and the constraint of \eqref{eq:ct_problem} by $c=1/(I_0\|W\|^2)$, which does not modify the set of
solutions, i.e., we take $g=c\|\cdot\|_{1,2}$ and $\tilde{d}=c\,d$, whose gradient is $1$-Lipschitz
continuous on $\RP^N$. If $\bar{w}\in\RPP$ is a Lagrange multiplier associated to the constraint of
this normalized problem, a solution $\bar{x}$ to \eqref{eq:ct_problem} is also a solution to
$\min_{x\in\RP^N}\lambda\,\mathrm{TV}(x)+\mathrm{KL}(b\,|\,\mu(x))$ with $\lambda=1/\bar{w}$, which is
the regularization parameter selected by the discrepancy principle.

We compare\footnote{The algorithms were implemented in Python~3 and run on
a computer equipped with an Intel Core i7-14700K processor (3.4/5.6~GHz), 64~GB of RAM, and running Windows~11 Pro 64-bit.}
Algorithm~\ref{algo:algononlinear} with the algorithms proposed in
\cite{BricenoDavis2018,BricenoRoldan4op}, which we denote by FBHF-LS and FBHF-LLS, respectively.
Both FBHF-LS and FBHF-LLS are primal algorithms applied to a primal-dual formulation of
\eqref{eq:ct_problem}, in which the linear composite term TV in \eqref{eq:ct_problem} is split by
means of the skew operator $\widetilde{D} \colon (x,u) \mapsto (\nabla^*u,-\nabla x)$. FBHF-LS includes
$\widetilde{D}$ in the continuous operator activated in the line search, while FBHF-LLS considers
$\widetilde{D}$ as an additional Lipschitz operator which is not involved in the line search
procedure. The advantage of Algorithm~\ref{algo:algononlinear} is that it is a primal-dual
algorithm, so it does not need to consider the linear term as a Lipschitz or continuous operator.

\begin{figure}[htbp]
\centering
\setlength{\tabcolsep}{1pt}
\begin{tabular}{ccccc}
\footnotesize Phantom $\overline x$ & \footnotesize Sinogram $\hat{p}$ & \footnotesize FBHF-LLS & \footnotesize FBHF-LS & \footnotesize Alg.~\ref{algo:algononlinear} \\
\includegraphics[width=0.19\textwidth]{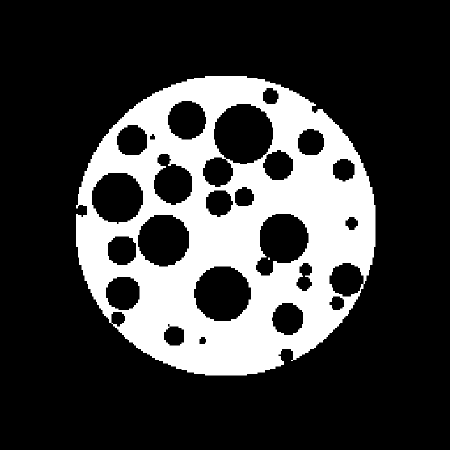} & \includegraphics[width=0.19\textwidth]{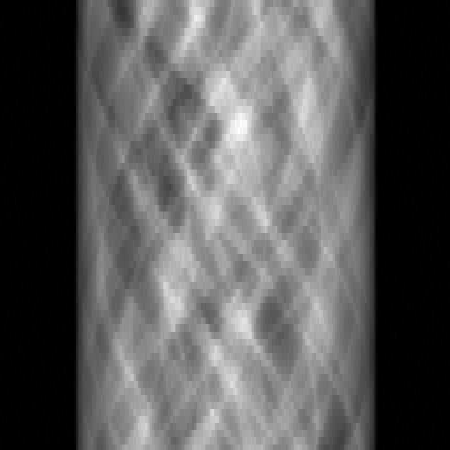} & \includegraphics[width=0.19\textwidth]{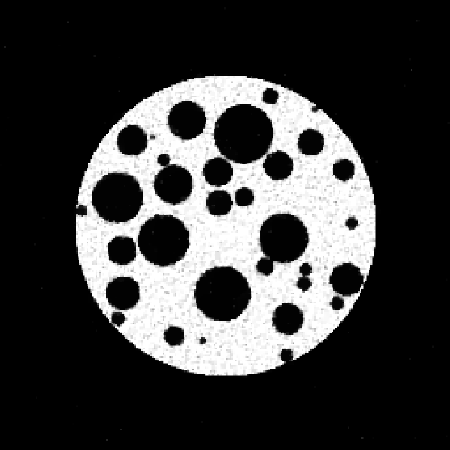} & \includegraphics[width=0.19\textwidth]{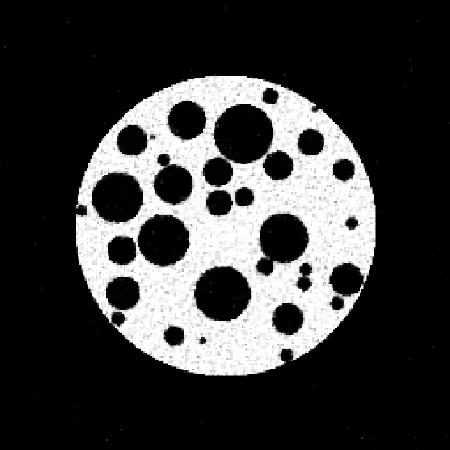} & \includegraphics[width=0.19\textwidth]{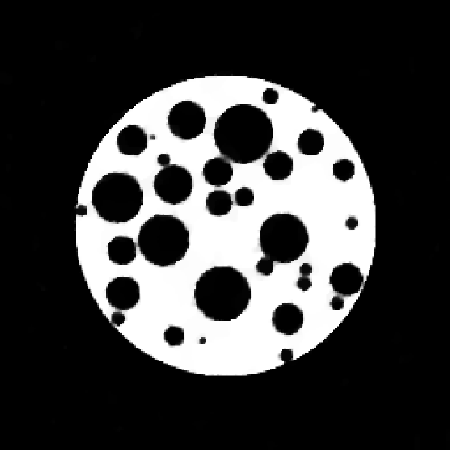} \\
\footnotesize $n=256$, $n_\theta=64$ & \footnotesize after 60 s & \footnotesize 24.21 dB, 9545 it. & \footnotesize 23.77 dB, 8615 it. & \footnotesize 26.87 dB, 9410 it. \\[2pt]
\includegraphics[width=0.19\textwidth]{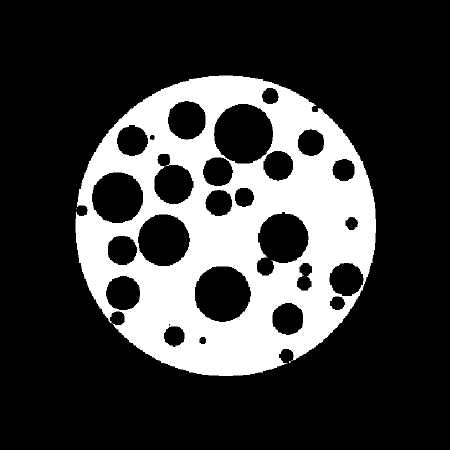} & \includegraphics[width=0.19\textwidth]{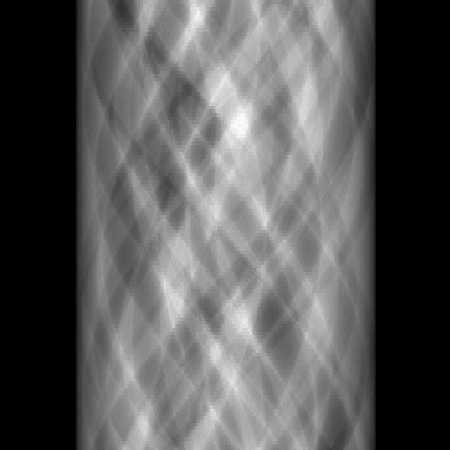} & \includegraphics[width=0.19\textwidth]{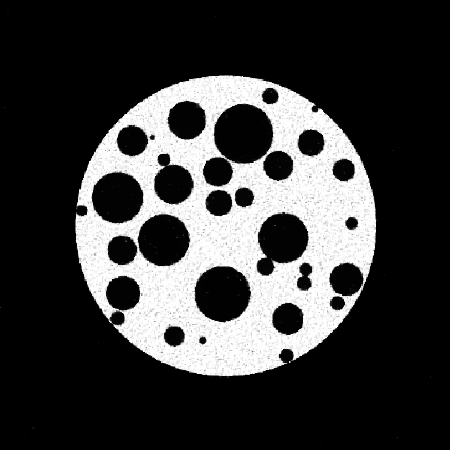} & \includegraphics[width=0.19\textwidth]{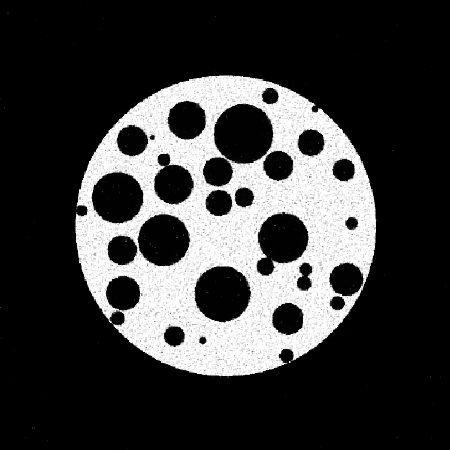} & \includegraphics[width=0.19\textwidth]{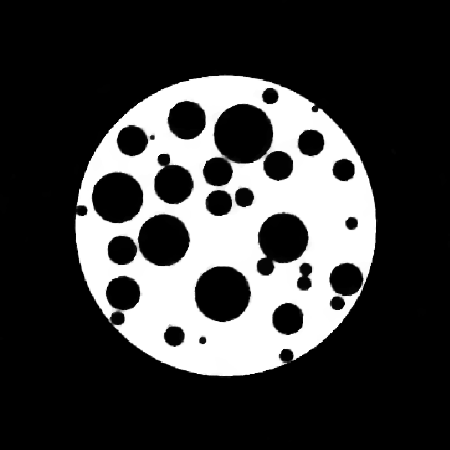} \\
\footnotesize $n=512$, $n_\theta=128$ & \footnotesize after 75 s & \footnotesize 24.94 dB, 11597 it. & \footnotesize 23.99 dB, 10246 it. & \footnotesize 28.40 dB, 10686 it. \\[2pt]
\end{tabular}
\caption{From left to right: phantom $x^\dagger$, measured sinogram $\hat{p}=\ln(I_0/\max\{b,1\})$, and reconstructions obtained by FBHF-LLS, FBHF-LS, and Algorithm~\ref{algo:algononlinear} ($\omega=0.01$), with warm line searches, after 60 seconds for $n=256$, $n_\theta=64$ and after 75 seconds for $n=512$, $n_\theta=128$. Below each reconstruction, the PSNR with respect to $x^\dagger$ and the number of iterations performed are shown.}
\label{fig:ct_images}
\end{figure}

We consider a two-dimensional foam-like phantom $\overline{x}$
generated with the package \texttt{foam\_ct\_phantom} \cite{Peltetal2022}, scaled such that
$\max_i(W\overline{x})_i=3$, and $I_0=10^4$. Hence, the expected counts satisfy
$\mu_i(\overline{x})\ge I_0e^{-3}\approx 498$, which corresponds to the large count regime considered in
\cite{Berteroetal2010}. The operator $W$ corresponds to a parallel-beam geometry with $n_\theta$
projection angles equally spaced in $[0,\pi[$ and $m=n$ pixels per side of the squared image and detector bins, so that $M=n n_\theta$. We
consider two scenarios: (i) $n=256$ and $n_\theta=64$, i.e., $N=65,536$ and $M=16,384$, and (ii)
$n=512$ and $n_\theta=128$, i.e., $N=262,144$ and $M=65,536$. In both cases, the data $b$ are
generated as independent Poisson random variables with means $\mu(\overline{x})$, and we set
$\eta=M/2$. All the methods start from the same point
$x_0=P_{\RP^N}\big(c_0W^\top\hat{p}\big)$, where $\hat{p}_i=\ln\big(I_0/\max\{b_i,1\}\big)$ is an
estimate of $(W\overline{x})_i$ obtained from the data, and $c_0\in\RPP$ is chosen such that
$\|c_0WW^\top\hat{p}-\hat{p}\|$ is minimal, that is, $x_0$ is a scaled (unfiltered) backprojection of
the log-data. The initial dual variables are $u_0=0$ and $w_0=1$ (for
Algorithm~\ref{algo:algononlinear}, we set $x_{-1}=z_0=x_0$ and $y_0=w_0$).
The phantom and the measured sinograms are shown in
Figure~\ref{fig:ct_images}.

We stop
the methods by using the residual of the primal-dual inclusion associated to
\eqref{eq:ct_problem}, i.e., for
Algorithm~\ref{algo:algononlinear}, we use $\rho_n=(s_n,t_n)$ defined in \eqref{eq:convtoZ}, and the method stops when 
$\|\rho_n\|\le10^{-6}$. A similar residual stopping criterion is used for FBHF-LLS and FBHF-LS. 
The operators $W$ and $W^\top$ were computed
on a GPU with the toolbox \texttt{tomosipo} \cite{Hendriksenetal2021} in single precision, and the
Kullback--Leibler divergence was evaluated in double precision.

We consider the following parameters. For Algorithm~\ref{algo:algononlinear} we set $\theta=0.707$, $\gamma=0.6$, and $\tau_{-1}=\tau=0.99\,\chi_{\omega}(\theta)=0.99\sqrt{(1-\theta^2)/(\omega\|\nabla\|^2)}$. 
Hence, a small value of $\omega$ allows for large primal
steps, but it also yields a large initial trial step-size. We report the results for
$\omega\in\{0.01,0.1\}$. 
For FBHF-LLS, since $\beta=+\infty$, we have
 $\rho=\sqrt{1-\varepsilon}/\|\nabla\|$, and we set $\varepsilon=0.1$, $\sigma=0.5$, and
 $\theta=(\sqrt{1-\varepsilon}-\|\nabla\|\rho\sigma)/2=\sqrt{1-\varepsilon}(1-\sigma)/2$ (see \cite[Theorem~3.3]{BricenoRoldan4op}). For FBHF-LS, since the
 cocoercive operator vanishes, $2\beta\varepsilon$ can be chosen freely and we set it as
 $1/\|\nabla\|$, $\varepsilon=0.1$, $\theta=\sqrt{1-\varepsilon}/2$, and $\sigma=0.6$ (see
 \cite[Theorem~2.3]{BricenoDavis2018}).  For a 2D image,  $\|\nabla\|$ can be estimated as $\|\nabla\|^2 \leq 8$ and $\|W\| $ was calculated with the power method.

For the three methods, the line search at iteration $n$
tests the step-sizes $\alpha_0,\alpha_0s,\alpha_0s^2,\ldots$, where $s\in\left]0,1\right[$ is the
backtracking factor ($s=\sigma$ for FBHF-LLS and FBHF-LS, and $s=\gamma$ for
Algorithm~\ref{algo:algononlinear}), and each trial step-size requires one forward and one back
projection. In the original line searches, which we call \emph{fixed}, the first trial step-size
is the same at every iteration, namely $\alpha_{\max}$.
Since $\alpha_{\max}$ can be much larger than the step-sizes
accepted in most of the iterations, we also consider a \emph{warm start} of the line searches,
similar to the strategy in \cite{MalitskyPock2018}, in which the first trial step-size, $\alpha_0$, is the last
accepted step-size enlarged by a factor $\delta>1$, that is, $
	\alpha_0=\min\{\delta\tau_{n-1},\alpha_{\max}\}$. Note that every trial step-size is bounded above by $\alpha_{\max}$, which is the only upper bound required for the convergence of the three methods. Moreover, a step-size is reduced only when the first trial fails. This strategy reduces the number of line search steps for each algorithm, thereby accelerating their convergence as we will see later. We set $\delta=1.5$.

\begin{figure}[htbp]
\centering
\includegraphics[scale=0.6]{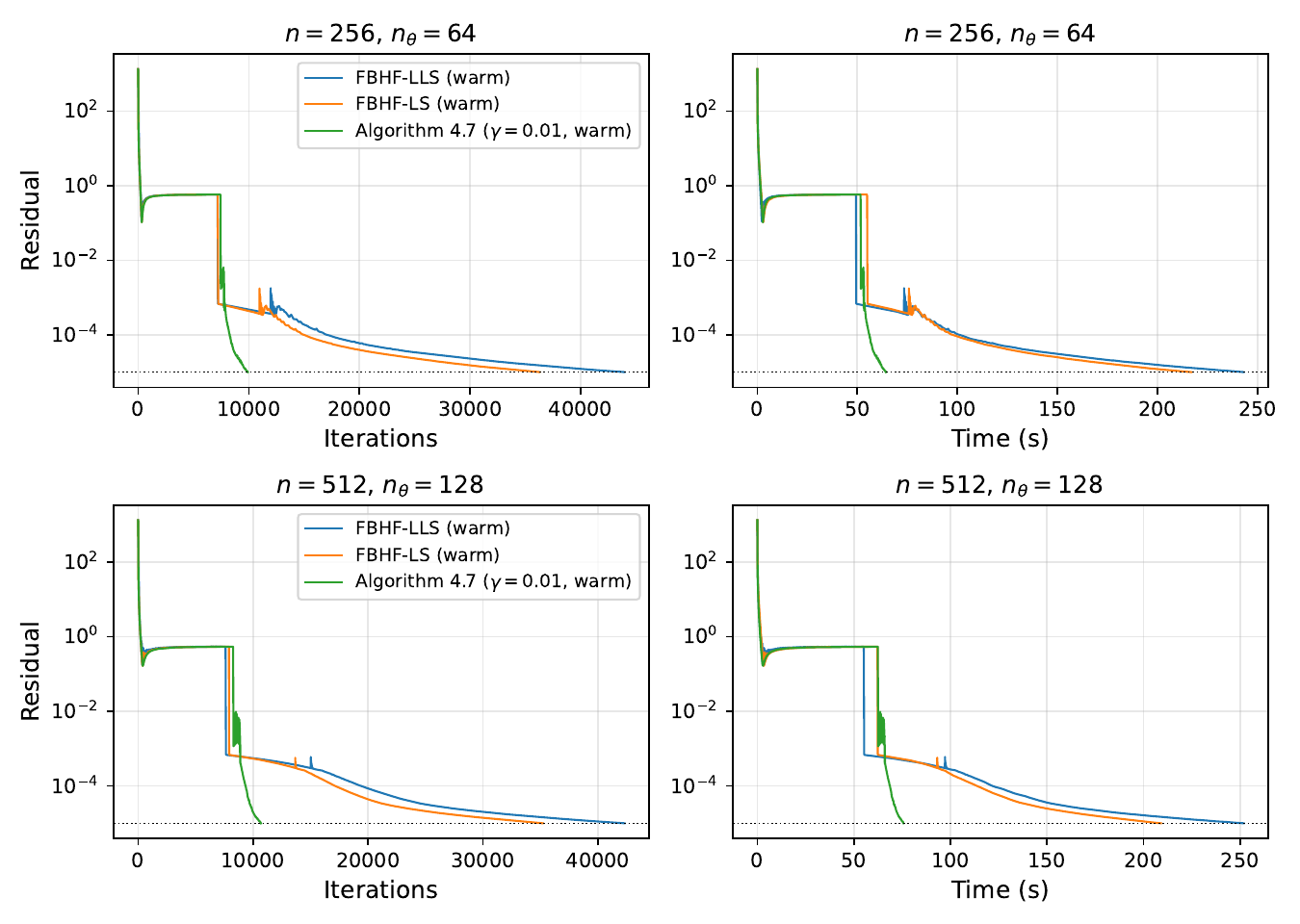}
\caption{Residual versus iterations (left) and time in seconds (right) for
FBHF-LLS, FBHF-LS, and Algorithm~\ref{algo:algononlinear}, with warm line searches, in scenarios (i) $n=256$, $n_\theta=64$ (top) and (ii)
$n=512$, $n_\theta=128$ (bottom). The dotted line is the tolerance $\epsilon=10^{-6}$.}
\label{fig:ct_residuals}
\end{figure}

Table~\ref{tab:results_ct} shows the results of the three methods with the
fixed and warm line searches, Figure~\ref{fig:ct_residuals} shows the evolution of the residuals for the fastest version of each method, and Figure~\ref{fig:ct_images} shows the
corresponding reconstructions. With the fixed line searches, Algorithm~\ref{algo:algononlinear}
needs 63--77\% fewer iterations than FBHF-LLS and 56--73\% fewer iterations than FBHF-LS, although its
advantage in computational time is smaller: in scenario (i), it is 19--26\% faster than FBHF-LLS and
23--30\% faster than FBHF-LS, while in scenario (ii), it is 14--15\% faster than FBHF-LS and it is
only 6--7\% faster than FBHF-LLS. The reason is that the
computational time is essentially proportional to the number of projections (column P), since the
products with $\nabla$ and $\nabla^*$ are much cheaper, and P is approximately IT$\times$(LS$+1$), where IT is the number of iterations and LS is the average number of line search trial step-sizes per iteration. Hence, the reduction in the number of iterations of Algorithm~\ref{algo:algononlinear} is
compensated by the larger number of trial step-sizes per iteration: 4.6--9.5, in comparison with
1.5--2.3 for FBHF-LLS and FBHF-LS, which is explained by the initial trial step-size
$\tau=0.99\,\chi_{\omega}(\theta)$, which is much larger than most of the accepted step-sizes.

\begin{figure}[htbp]
\centering
\includegraphics[scale=0.5]{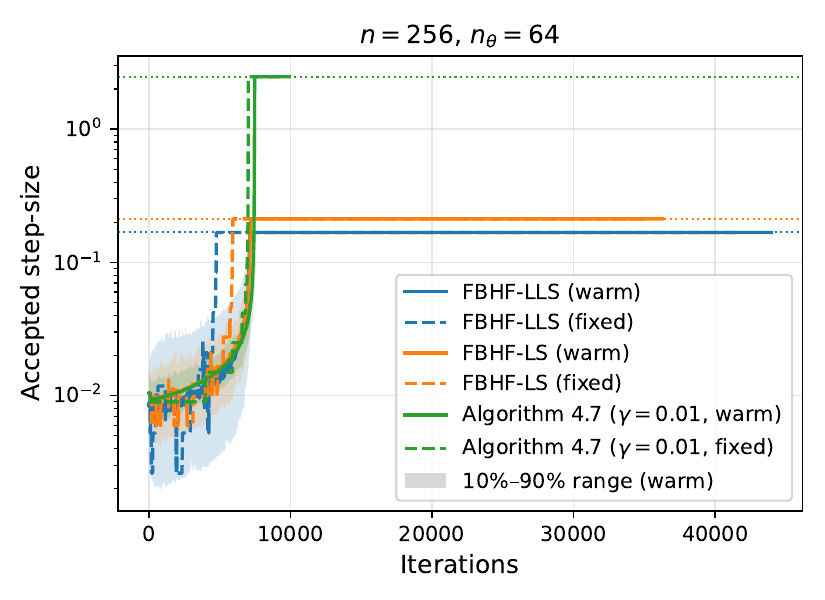}\ 
\includegraphics[scale=0.5]{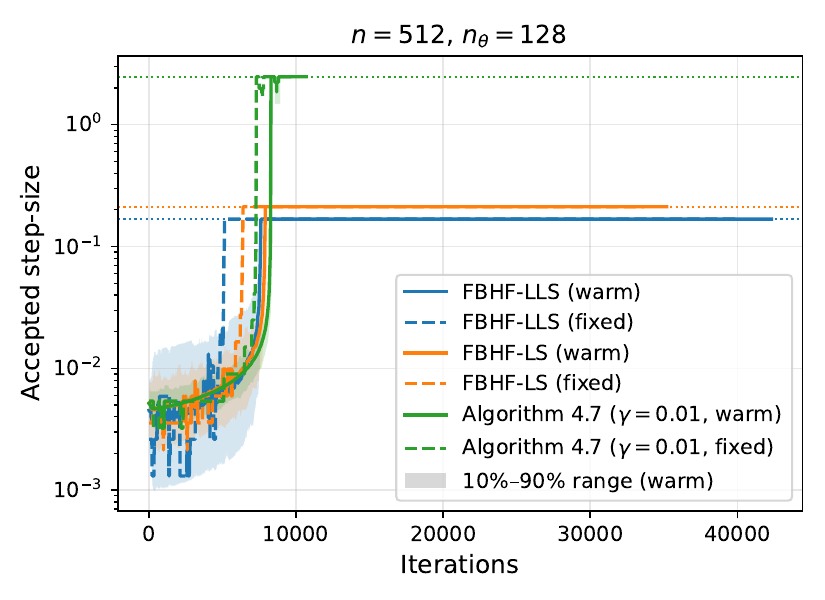}
\caption{Accepted step-sizes versus iterations (median over a moving window of 200 iterations) for
FBHF-LLS, FBHF-LS, and Algorithm~\ref{algo:algononlinear}. The dotted lines are the maximal trial step-sizes $\alpha_{\max}$
of each method.}
\label{fig:ct_steps}
\end{figure}

The warm start solves this issue. For Algorithm~\ref{algo:algononlinear} with
$\omega=0.01$, it reduces the average number of trial step-size per iteration from 8.4--9.5 to
1.6, the number of projections by 71--73\%, and the computational time by 70--73\%, while the number
of iterations increases by only 5--10\%. For FBHF-LLS and FBHF-LS, since $\alpha_{\max}$ is close
to the accepted step-sizes, the fixed line searches already need few trial step-sizes and the warm
start reduces the computational time by 8--15\% and 23--35\%, respectively. Note that the warm start can
hardly reduce the number of iterations of FBHF-LLS and FBHF-LS, since their step-sizes are limited by
$\|\nabla\|$, because the linear operator is treated as a Lipschitz operator. Indeed,
their numbers of iterations increase by 3--6\%. Hence, with the warm
start, Algorithm~\ref{algo:algononlinear} with $\omega=0.01$ is the fastest method in both
scenarios: it needs 70--78\% fewer iterations than FBHF-LLS and FBHF-LS, and it is 3.3--3.8 and
2.8--3.4 times faster than FBHF-LLS and FBHF-LS, respectively. Moreover, with the warm start, the
smallest value of $\omega$ provides the best performance, since the large initial trial step-size
$\chi_\omega(\theta)$ no longer leads to wasted trial step-sizes, which simplifies the choice of
$\omega$. Figure~\ref{fig:ct_residuals} shows that  Algorithm~\ref{algo:algononlinear} reaches the tolerance in fewer iterations and less GPU time. This advantage can be explained by the fact that Algorithm~\ref{algo:algononlinear} allows for larger step-sizes.
Figure~\ref{fig:ct_steps} plots the moving median of the accepted step-sizes over 200 iterations. We observe that in the initial phase of the iterations, all algorithms iterate with a similar step-size, which increases until reaching their maximum trial step-sizes. Indeed, this figure shows that the largest trial step-size is provided by Algorithm~\ref{algo:algononlinear}, thereby accelerating convergence. Figure~\ref{fig:ct_images} shows the reconstructions obtained by the three methods after 60 seconds for $n=256$ and 75 seconds for $n=512$. From this figure, we observe that Algorithm~\ref{algo:algononlinear} provides a better reconstruction within these execution times.
{
\begin{table}[htbp]
\centering
\resizebox{13cm}{!}{
\begin{tabular}{cccccccccc}
\toprule
 &  & \multicolumn{4}{c}{$n=256$, $n_\theta=64$} & \multicolumn{4}{c}{$n=512$, $n_\theta=128$} \\
\cmidrule(lr){3-6} \cmidrule(lr){7-10}
Line Search & Algorithm & T & IT & P & LS & T & IT & P & LS \\
\midrule
Fixed & FBHF-LLS & 265.04 & 41,680 & 102,154 & 1.45 & 295.70 & 39,861 & 105,001 & 1.63 \\
 & FBHF-LS & 281.00 & 35,091 & 102,005 & 1.91 & 323.32 & 33,661 & 110,135 & 2.27 \\
 & Alg.~\ref{algo:algononlinear} ($\omega=0.01$) & 215.90 & \textbf{9,477} & 89,223 & 8.41 & 278.35 & \textbf{9,754} & 102,539 & 9.51 \\
 & Alg.~\ref{algo:algononlinear} ($\omega=0.1$) & 197.47 & 14,311 & 79,795 & 4.58 & 275.01 & 14,942 & 100,380 & 5.72 \\
\midrule
Warm ($\delta=1.5$) & FBHF-LLS & 243.08 & 43,965 & 92,086 & 1.09 & 251.90 & 42,312 & 89,065 & 1.10 \\
 & FBHF-LS & 216.95 & 36,291 & 78,261 & 1.16 & 209.04 & 35,204 & 76,668 & 1.18 \\
 & Alg.~\ref{algo:algononlinear} ($\omega=0.01$) & \textbf{64.28} & 9,909 & 25,733 & 1.60 & \textbf{75.70} & 10,686 & 28,064 & 1.63 \\
 & Alg.~\ref{algo:algononlinear} ($\omega=0.1$) & 93.23 & 15,252 & 36,488 & 1.39 & 95.10 & 15,382 & 37,228 & 1.42 \\
\bottomrule
\end{tabular}}
\caption{Time in seconds (T), number of iterations (IT), number of forward projections (P), 
and average number of trial step-sizes per iteration (LS), needed to
satisfy the stopping criteria, for FBHF-LLS,
FBHF-LS, and Algorithm~\ref{algo:algononlinear} with $\omega\in\{0.01,0.1\}$, fixed and warm line searches, in scenarios (i) $n=256$ with $n_\theta=64$
projection angles and (ii) $n=512$ with $n_\theta=128$ projection angles.}
\label{tab:results_ct}
\end{table}}

\section{Conclusion}\label{sec:conclu}
In this article, we proposed a splitting algorithm with a primal--dual line search for solving monotone
inclusions involving a maximally monotone operator, a cocoercive operator, a continuous operator, and
a composite term with a bounded linear operator. To the best of our knowledge, this is the first line
search in which the variation of a merely continuous operator is controlled by a joint primal-dual
displacement. This test, which evaluates only the continuous operator, allows us to handle the linear
operator as in the Condat--V\~u and Chambolle--Pock methods. The step-sizes adapt to the local behavior
of $D$, and no Lipschitz constant of $D$ is needed. We showed by examples that the dual term in this test cannot be dropped, neither for the termination
of the line search nor for the convergence of the iterates.
We proved the weak convergence of the iterates under mild assumptions on the continuous
operator, recovering several known methods in the Lipschitz case, and we applied our results to convex
problems involving $K$-convex nonlinear compositions. In computed tomography experiments, the proposed method needed
considerably fewer iterations than FBHF-type methods with line search, since it allows for larger
step-sizes, and, combined with a warm-started line search, it was significantly faster. 

\section*{Acknowledgment}
The work of Felipe Atenas and Luis M. Brice\~no-Arias is supported by  Centro de Modelamiento Matem\'atico (CMM), ACE210010
and FB210005, BASAL fund for 
centers of excellence. Luis M. Brice\~no-Arias was also supported by FONDECYT 1230257 from ANID-Chile. The work of Fernando Roldán is supported by ANID-Chile under grant FONDECYT Iniciación 11250164.

\bibliographystyle{spmpsci}
\bibliography{ref}

\end{document}